# Equity-aware Design and Timing of Fare-free Transit Zoning under Demand Uncertainty

Qianwen Guo
Department of Civil and Environmental Engineering,
Florida State University
2525 Pottsdamer St, Tallahassee, FL 32310
Email: qguo@fsu.edu

Jiaqing Lu
Department of Civil and Environmental Engineering,
Florida State University
2525 Pottsdamer St, Tallahassee, FL 32310
Email: jl23br@fsu.edu

Joseph Y. J. Chow
Department of Civil and Urban Engineering,
New York University Tandon School of Engineering
6 MetroTech Center Brooklyn, NY 11201, United States
Email: joseph.chow@nyu.edu

Paul Schonfeld (corresponding author)
Department of Civil and Environmental Engineering,
University of Maryland
1173 Glenn Martin Hall, College Park, MD 20742, United States
Email: pschon@umd.edu

# Abstract

We propose the first analytical stochastic model for optimizing the configuration and implementation policies of fare-free transit. The model focuses on a transportation corridor with two transportation modes: automobiles and buses. The corridor is divided into two sections, an inner one with fare-free transit service and an outer one with fare-based transit service. To capture the uncertainty in transit demand, the model uses a Geometric Brownian Motion process to represent the evolution of demand density, thus enhancing the realism of the analysis. Under the static version of the model, the optimized length and frequency of the fare-free transit zone can be determined by maximizing total social welfare. The findings indicate that implementing fare-free transit can increase transit ridership and reduce automobile use within the fare-free zone while social equity among the demand groups can be enhanced by lengthening the fare-free zone. Notably, the optimal zone length increases when both social welfare and equity are considered jointly, compared to only prioritizing social welfare. The dynamic model, framed within a market entry and exit real options approach, solves the fare policy switching problem, establishing optimal timing policies for activating or terminating fare-free service. Closed-form analytical solutions are derived for the implementation policies, including a switching cost-based-threshold (also called timing) for fare-free transit service, a switching cost-based-threshold for fare-based service, and a single threshold without switching costs under social welfare- and equity-aware regimes, respectively. We applied the dynamic models to Utah Transit Authority (UTA) Blue Line to optimize the timing for implementing fare-free transit and the length of fare-free transit zone under social equity-aware and social welfare-aware regimes. The results from dynamic models reveal earlier implementation and extended durations of fare-free transit in the social welfare-aware regime, driven by lower thresholds compared to the social equity-aware regime. The proposed tool provides valuable insights for policymakers and transit planners, enabling them to implement fare-free transit services more effectively from social welfare and equity perspectives.

**Key Words***:* Fare-free transit, Equity, Implementation policies, Optimal design, Demand uncertainty.

# 1. Introduction

The COVID-19 pandemic has stimulated conversations around the potential benefits and impacts of fare-free transit. Many transit agencies eliminated fares during the pandemic, with the initial intention of reducing contact between riders and operators (Dai et al., 2021; Goldberg, 2021). As life largely returned to normal, some transit agencies, such as DASH in Alexandria, Virginia, ABQ Ride in Albuquerque, New Mexico, or the entire country of Luxembourg, extended free fares in the post-pandemic era, hoping to recover lost ridership and financially support transit riders in the wake of the economic downturn (Kirschen et al., 2022).

Transit agencies in the U.S. and elsewhere have long experimented with fare-free transit services due to their potential advantages related to cost, system efficiency, access, and operations (Goldberg, 2021). When fares are waived, the whole fare collection system can be discarded, along with its associated costs, such as for hardware and staffing to manage and count fares. For instance, Intercity Transit in Olympia, WA found that eliminating fareboxes could avoid long-run costs of system upgrades. Fares should not be collected when the collection costs exceed the resulting revenues, especially by transit agencies with relatively low farebox recovery ratios (Spegman, 2019). A literature scan shows that fare-free transit programs seem more successful in smaller transit systems, where the collection costs reduce net revenues disproportionately for shorter trips. This is demonstrated by long-running fare-free transit services in small cities (e.g., Commerce, California), resort communities (e.g., Breckenridge, Colorado), and university-dominated communities (e.g., Amherst, Massachusetts) (Volinski, 2012).

Determining whether fare-free transit service is adopted involves a delicate balance between benefits and costs. Fare-free transit has several benefits, as it increases household utility by reducing transportation expenses and eliminates costs associated with farebox collection (Chen & Zhou, 2022). It may also result in lost fare revenue, increased cost of operations, and reduced level of service due to increased ridership (Cantillo et al., 2022). For instance, Intercity Transit in Olympia, USA, found that eliminating fareboxes could save costs in the long run as no costly fare collection system upgrades would be needed. On the other hand, transit researchers have indicated that a fare-free policy may cause loss of fare revenue, increased cost of operations, reduction in the level of service, and potential security issues (Saphores, et al., 2020).

For a transit agency, it is challenging to decide whether to adopt a fare-free policy, how to implement it, when to implement it, and when to terminate it, because such decisions depend on the range of services a transit agency provides, travel patterns, socioeconomics of its passengers, and financial support that agency receives from its sponsor (Lu et al., 2024). At one extreme, none of the transit services are free; at the other extreme, all services are free. Between those two extremes, there is a broad range of alternatives. Transit operators can decide whether to introduce a fare-free policy and fine-tune their policy by considering service types, areas of coverage, operating periods, and demographics (such as income). It should be noted that those policy schemes are not mutually exclusive and can be bundled, depending on the needs and objectives of a transit agency.

The diversity of policy schemes is also consistent with practice. For smaller transit systems, the policies can be adopted system wide. It is common for larger transit agencies to offer fare-free transit partially, which means they waive fares only for specific groups of riders, on certain routes or transit modes, during certain time periods, or in pre-defined areas. For example, the Utah Transit Authority has adopted a partial fare-free structure which allows bus passengers to ride for free only in dedicated zones (Regional Zero-Fare Study, 2022). In other instances, a few municipalities in Florida provide fare-free circulator services within their jurisdictions, such as Metromover in Miami-Dade County, the Bus Rapid Transit (BRT) circulator LYMMO in Orange County, and the beach trolley and tram service provided by LeeTran in Lee County (Kębłowski, 2020).

The decision to initiate or discontinue fare-free transit services typically hinges on a social welfare analysis conducted by transit agencies (Butler & Sweet, 2020). The feasibility of such a policy is significantly influenced by the density of the demand it serves. Nevertheless, the future demand density along a transportation corridor is often uncertain, a fact underscored by the COVID-19 pandemic and its significant impact on transit ridership. In response to decreased ridership, some agencies have elected to discontinue fare-free transit, while others have sustained fare-free transit services to provide economic benefits to patrons. The deliberation on whether to launch or terminate fare-free transit, the timing of these actions, and the service operation decisions, remain critical considerations for many transit agencies.

In this study, the implementation policies grounded in the market entry-exit real options approach is investigated, relying on a stochastic demand model. Dixit (1989) initially developed this approach to resolve certain financial investment problems. Subsequently, Sødal et al. (2008) adapted this approach for maritime shipping. The bulk shipping market segregates into two sectors, determined by cargo type, namely dry bulk and wet bulk. A specific type of ship, called a combination carrier, can transport both types of cargo. The owner of a combination carrier determines when to transition from wet bulk operations to dry bulk operations, considering a given switching cost and based on the current freight rates in both markets. Guo et al. (2018) further applied this method to optimize various transit service switching problems and to handle issues related to operating shared autonomous transit fleets.

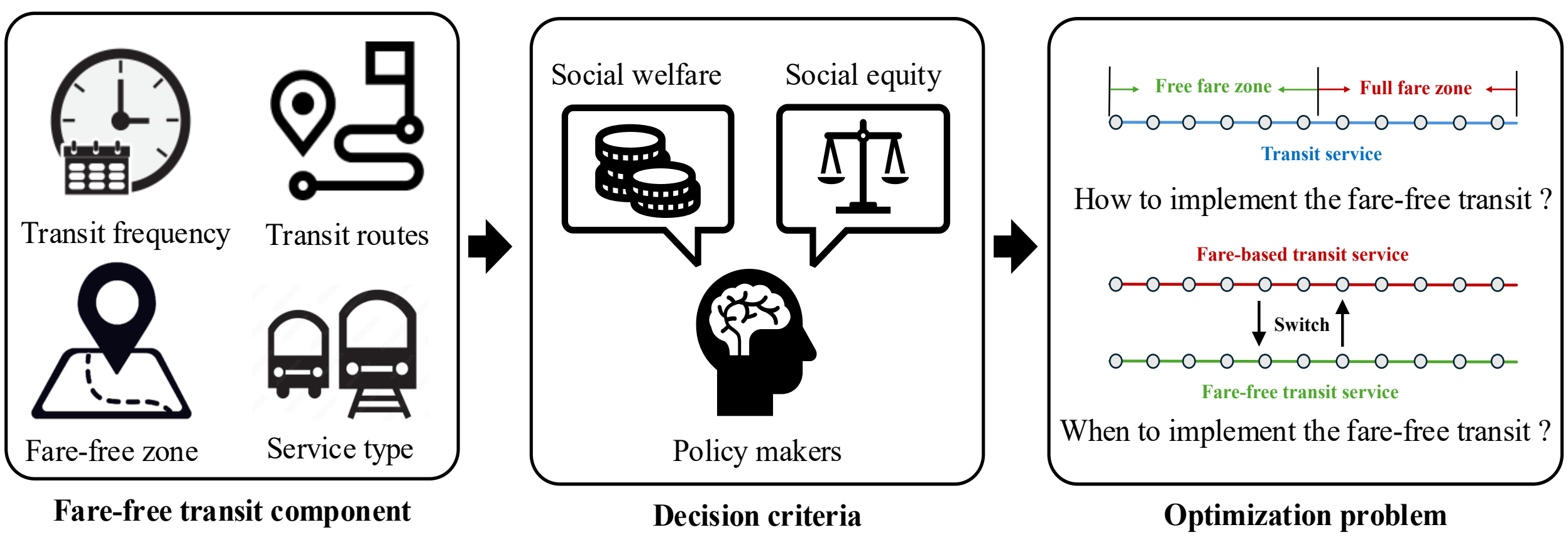

Figure 1 Optimized design of fare-free transit

A survey of existing literature on fare-free transit reveals significant gaps in quantitative analyses and optimization models pertaining to the spatial design of transit systems. Furthermore, there is a lack of dynamic models that assess the feasible periods for fare-free transit over an extended planning period. Additionally, the current body of research on optimal timing offers limited insights into the simultaneous determination of fare-free transit zone allocation and implementation strategies.

This study seeks to address these deficiencies and make the following contributions:

a. It aims to optimize the length of fare-free transit service zones and the service frequency by maximizing total social welfare across static, dynamic, and stochastic dynamic frameworks. The latter two frameworks incorporate mechanisms for activating and deactivating fare-free operations, leveraging an original real options market entry-exit switching policy. The models consider scenarios both with and without social equity consideration, focusing on the context of a bimodal transportation system. Analysis of such systems reveals insights on discontinuities in mode split and user surplus due to the zone design.

b. It quantitatively assesses the social equity implications of implementing fare-free transit services. It examines the disparities among residents along a commuting corridor, by analyzing the impact of fare-free transit service on the user surplus of demand groups along the corridor. A benefit index is proposed to balance social welfare and social equity, and the length of the free zone is optimized to maximize overall benefit. Analysis reveals a nonlinearity in the dependency of the benefit index on the design.
c. It tests the stochastic policy by simulating dynamic stochastic demand with data from the UTA Blue Line. The results verify the effectiveness of this policy.

The justification for this paper lies in its potential contributions to the field. By exploring the intricate relations among transit system operation, user cost parameters, and the spatial design of fare-free transit zones, this research offers technical guidance for decision-makers. It aims to support more strategic and effective implementation of fare-free transit services than short-sighted or simplistic decisions. The findings from this study can inform policymakers and transportation authorities, enabling them to make good decisions regarding the implementation and optimization of fare-free transit systems. Ultimately, this research aims to enhance the accessibility, effectiveness, and equity of public transportation systems.

# 2. Literature Review

## 2.1 Fare-free transit examples and developments

The concept of fare-free transit, also referred to as zero-fare transit, is an innovative approach that breaks with conventional transportation funding policies by dispensing with passenger fares (Studenmund & Connor, 1982). The development timeline of fare-free transit is shown in Figure 2. The advent of fare-free transit can be traced back to at least 1960, when Hasselt, a Belgian city, embarked on an experiment in this arena (Saphores et al., 2020). In the 1970s, several universities began implementing free bus services for their students and employees, including UC Berkeley and the University of Maryland (Kirschen et al., 2022). Over the next two decades, several small municipalities worldwide, including Chapel Hill in North Carolina, USA, initiated tests of fare-free transit in their respective regions (Bull et al., 2021).

A surge in urban areas and regions considering or implementing fare-free transit has been observed in the 21st century. One notable instance is Tallinn, the capital of Estonia, where public transportation became fare-free for all registered inhabitants (Gray, 2018). In 2018, Dunkirk, France, also made a significant move by providing its buses and trams to locals and tourists free of charge. Compared with some large countries, a small country may more easily launch fare-free transit throughout its territory. Thus, Luxembourg marked a milestone in 2020 by becoming the first country to implement fare-free transit nationally (Saphores et al., 2020).

However, the onset of the COVID-19 pandemic has profoundly influenced the progression of fare-free transit. Due to the reduction in ridership triggered by the pandemic, many agencies have had to withdraw their fare-free schemes. Conversely, some transit authorities, in an attempt to maintain social distancing on vehicles and assure the mobility of essential workers, have chosen to offer free transit (Goldberg, 2021, Dai et al., 2021). In the post-COVID-19 world, the decision to uphold or abandon fare-free transit remains a topic of ongoing debate among policymakers, transportation experts, and communities.

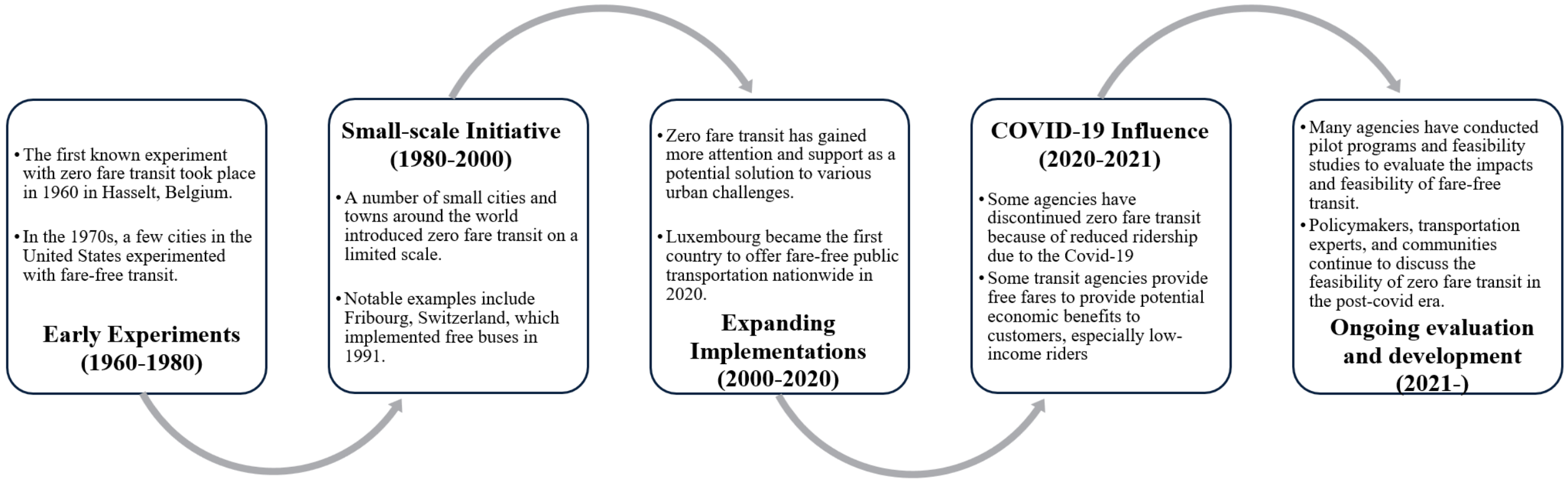

Figure 2: Timeline of fare-free transit development.

The adoption of fare-free transit introduces multifaceted impacts, particularly on ridership growth and shifts in transportation mode choices. De Witte et al. (2008) detected a definite increase in public transport usage following the introduction of free services; however, this must be accompanied by enhancements in the transport system's quality, such as its frequency, capacity, and connections. Boyd et al. (2003) scrutinized the modal shift at the University of California's Los Angeles campus after rendering bus transportation free of charge. The study found that transit ridership soared by over 50%, resulting in a reduction of more than 1000 daily car trips to the campus.

Similarly, De Witte et al. (2006) explored the implications of fare-free public transportation for students in Brussels, and observed an uptick in public transportation use when services were made costless. Verheyen (2010) assessed the influence of free public transport on modal splits, distinguishing between different trip motives, such as commuting, shopping, and recreational activities. Cats et al. (2017) evaluated the changes observed prior to and approximately a year after implementing fare-free transit using an annual municipal survey. Their findings indicated that the market share of public transportation climbed by 14%, attributable to a 10% decrease in car trips and a 40% drop in walking as the primary travel mode.

From these studies, it is clear that fare-free transit typically enhances ridership, presenting potential challenges to existing service schedules. Numerous academics have deliberated the possible drawbacks and costs associated with fare-free transit implementation on service schedules. Such schedules may need revisions to accommodate operational changes that could arise from zero- or reduced-fare options, for instance, expedited travel times due to decreased boarding durations or decreased service level resulting from increased ridership (Cats et al., 2017).

The introduction of fare-free transit also brings along a range of benefits, notably in terms of promoting social equity. Proponents of this approach argue that it could significantly enhance accessibility for vulnerable demographics and foster social equity (Cats et al., 2012). Social equity for different transit users is an important consideration when optimizing transit fares and subsidies (Wang et al., 2019; Wang et al., 2023). They designed a differential fare scheme which considers different groups in an urban rail transit system. To leverage the benefits of fare-free transit in promoting social equity more effectively, it would be useful to develop tailored fare-free policies to foster different equity types. For instance, Tallinn has set a clear objective of facilitating job search for its residents by offering accessible transportation (Cats et al., 2012). This policy can also improve the mobility of students, enabling their participation in cultural, social, and educational activities, as well as enhance the mobility of senior citizens, generating notable economic benefits for the local economy (Cats et al., 2012). Some theorists suggest that the act of purchasing a ticket, despite substantial discounts, can significantly dampen the propensity to utilize public transportation due to the psychological burden it imposes (Hodge et al., 1994; Noël et al., 2022).

## 2.2 Real options approach in transportation applications

The Real Options Approach (ROA) is a flexible decision-making tool that acknowledges the inherent uncertainties and irreversibility of investments. Its application in the transportation sector has significantly increased over the years. Indeed, ever since Dixit and Pindyck's seminal work on investment under uncertainty (D&P ref, 1994), the real options approach has been widely utilized to address a plethora of investment challenges in transportation infrastructure. The method solves an optimal timing problem under non-stationary uncertainty, where the stochastic process(es) varies over time, i.e., is characterized as a type of Markov decision process (Chow, 2018).

For complex multidimensional decisions (e.g., network design (Chow and Regan, 2011a, 2011b)), computational methods are used to find the optimal policy and value (e.g., Brennan and Schwartz, 1978; Cox et al., 1979; Longstaff and Schwartz, 2001), as a type of value function approximation (VFA) in approximate dynamic programming (Powell, 2007; Chow and Sayarshad, 2016). These methods typically employ a geometric Brownian motion (GBM) or mean-reverting process to model demand uncertainty. More recent work in this area includes using the VFA policy to monitor traffic over time (Chow, 2016) or integrating deep learning to design and time the expansion of service regions over time (Rath and Chow, 2022).

Alternatively, when the real options problem is sufficiently simple (single decision, single stochastic process) within an infinite horizon, optimal thresholds for making the decisions can be derived analytically either in closed form or through numerically-solved systems of equations. In addition to the deferral option for an irreversible decision presented in Dixit and Pindyck (1994), there is also a more generalized "switching option" in which decisions are made to switch reversibly (but with asymmetric switching costs) between two operating modes (Dixit, 1989). In the latter, with asymmetric switching costs, the presence of hysteresis extends one switching threshold to two.

Couto et al. (2015) applied the real option approach to a High-Speed Rail (HSR) project, deriving the demand threshold for its implementation. Guo et al. (2018) derived such a switching option for bi-operational public transit services considering mean reverting demand process. Balliauw et al. (2020) optimized the capacity investment decision for a privately-owned airport under demand uncertainty, modeling the aircraft movements with a GBM and selecting profit maximization as the optimization objective due to the airport's private status.

Guo et al. (2023) utilized the jump diffusion model to represent demand uncertainty in rail transit, jointly determining the timing and sizing of investments. However, it should be noted that previous research typically focused on timing issues, often overlooking other critical factors such as transit scheduling decisions (e.g., transit coverage length, service frequency), and the social equity among users distributed in different locations.

## 2.3 Summary

The empirical research on fare-free public transit is substantial, albeit with limited theoretical studies focusing on the optimization of transit pricing. Guo et al. (2021) established the viability of fare-free solutions that can improve social welfare, especially in instances where transit resource constraints are non-binding. However, a review of pertinent studies on fare-free public transit operations reveals several gaps. First, there is a lack of research on the optimization of fare-free transit zones and service schedules. Second, the thresholds for activating or discontinuing fare-free transit have not been adequately considered, particularly in more realistic scenarios where future demand remains uncertain. Finally, the issue of social equity for commuters along a transportation corridor has been largely overlooked when formulating fare-free transit service schemes and implementation policies.

This study aims to address these limitations by enhancing existing transit pricing research. It proposes optimized fare-free transit design and implementation policies that consider social welfare and social equity while incorporating stochastic demand evolution in the real-world case. The goal is to contribute to a more comprehensive understanding of fare-free transit, providing valuable insights for policy makers and transit planners.

# 3. Model Setup

To clearly define our model, all notations used in this research are listed in Table 1.

Table 1: Notation list

| Notation | Definition | Unit |
| --- | --- | --- |
| $A$ | City boundary | mile |
| $B$ | Fare-free transit zone length | mile |
| $C_a$ | Auto user travel cost | $ / trip |
| $C_b$ | Bus user travel cost | $ / trip |
| $C_f^a$ | Fixed cost of auto user | $ |
| $C_v^a$ | Variable cost for auto user | $/mile |
| $D$ | Cost of switching from fare-free transit to fare-based transit service | $/time |
| $E$ | Equity ratio of fare-free transit schemes | - |
| $e_c$ | Coefficient of generalized cost of auto users | - |
| $e_0$ | Fixed cost associated with the establishment of the fare collection system | $ |
| $e_1$ | Variable cost per unit of length of fare collection system | $/mile |
| $e_2$ | Fare collection cost per trip | $/trip |
| $F_i$ | Transit frequency at stage *i* | vehicle /hour |
| $f_i$ | Transit fare at stage *i* | $ |
| $G$ | Benefit index of fare-free transit schemes | - |
| $g_f$ | Fixed cost for transit operati0n | $/day |
| $g_v$ | Variable cost for transit operation | $/vehicle/day |
| $i$ | Transit service stage | - |
| $K$ | Cost of switching from fare-based transit to fare-free transit service | $/time |
| $k$ | Monthly discount rate | % |
| $N_i$ | Fleet size | vehicle |
| $n$ | The number of demand groups spatially distributed along the corridor | - |
| $P$ | Demand groups spatially distributed along the corridor extended from the CBD to the city boundary | pax/mile/day |
| $Pr_b$ | Probability of selecting bus by users | % |
| $Pr_a$ | Probability of selecting auto by users | % |
| $Q_{CBD}$ | Current demand density at central business district (CBD) | pax/mile/day |
| $\overline{Q}$ | Upper bound demand threshold from no fares to fares | pax/mile/day |
| $\underline{Q}$ | Lower bound demand threshold from fares to no fares | pax/mile/day |
| $R_i$ | Fare revenue at stage *i* | $ |
| $S$ | Average access time to transit stops | hour |
| $T_i$ | Expected discounted social welfare without ($i = 0$) and with ($i = 1$) fare-free | $ |
| $U_i$ | User surplus | $ |
| $\overline{U_{ab}}^{\zeta}$ | The average user surplus of demand group $\zeta$ | $ |
| $V_0(Q)$ | Option value of operating fare-based transit service | - |
| $V_1(Q)$ | Option value of implementing fare-frees service | - |
| $W_i$ | Total social welfare per day at stage *i* | $ |
| $W_{fi}$ | Fare collection cost at stage *i* | $ |
| $W_{oi}$ | Operating costs per day | $ |
| $W_{d1}$ | Administrative costs of fare-free program | $ |
| $\alpha_T$ | Value of in-vehicle travel time | $/hour |
| $\alpha_w$ | Value of waiting time | $/hour |
| $\beta_{min}, \beta_{max}$ | Social equity constraints | - |
| $\iota_f$ | Fixed administration cost of fare-free program | $ |
| $\iota_v$ | Variable administration cost of administration cost for fare-free program | $/mile |
| $\theta$ | Parameters to reflect the nonlinear nature of the cost growth | - |
| $\eta$ | Monthly growth rate of CBD demand | % |
| $\sigma$ | Monthly volatility rate CBD demand | % |
| $\psi$ | The parameter in the Logit model | - |

## 3.1 Assumptions

To facilitate the presentation of the essential ideas without losing generality, the following basic assumptions are made.

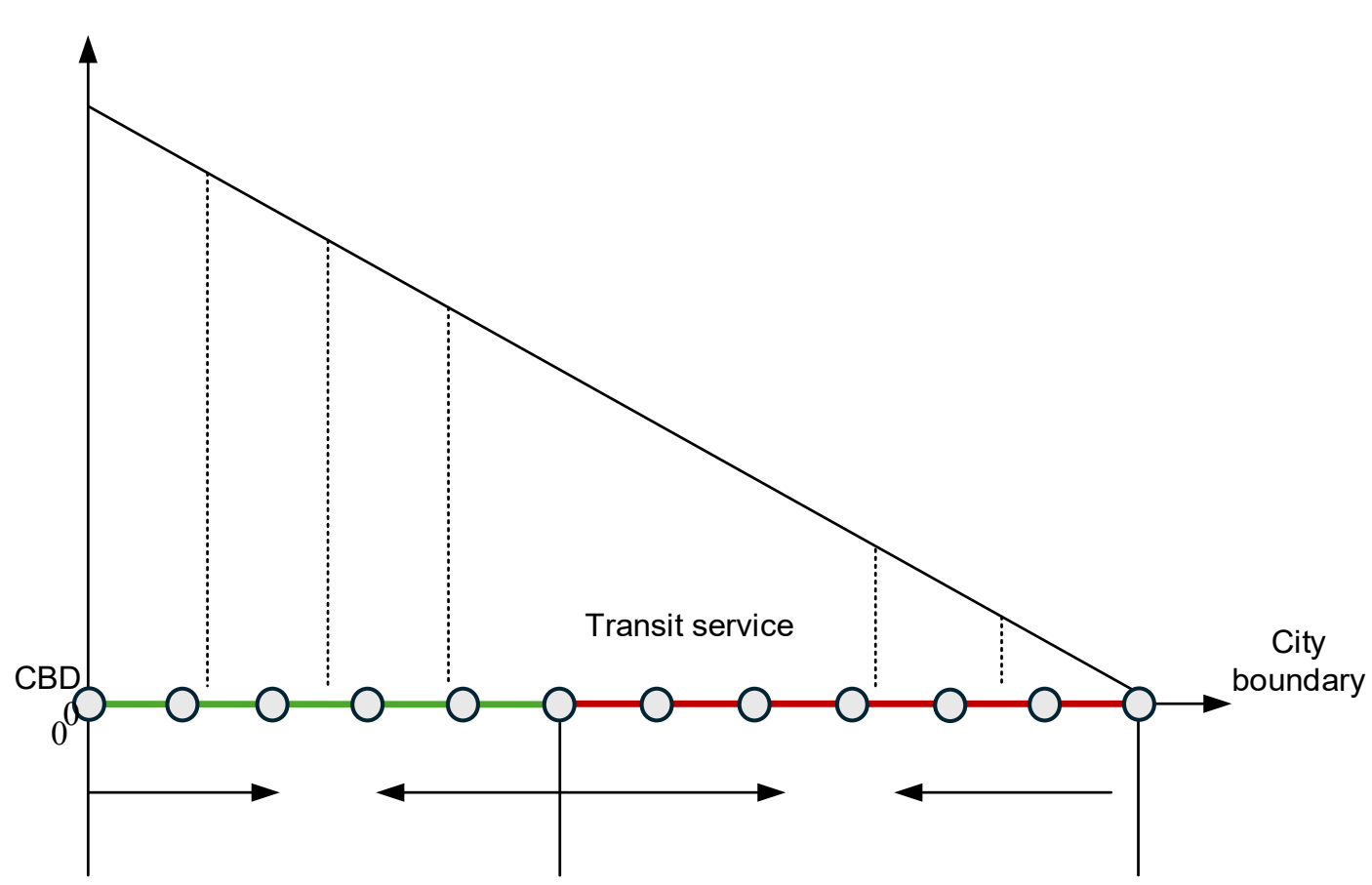

Figure 3: Demand density along the transportation corridor.

**A1** The corridor connecting the city's CBD and suburban area is assumed to be linear, as shown in Figure 3. Two alternative modes, auto and bus, share the same roadway along this transportation corridor. The bus line operates along the transportation corridor, which is subdivided into two sections: From the CBD the inner section $B$ is covered by fare-free transit service ($F_F$-service), while the outer section $A - B$ to the city boundary $A$ is covered by fare-paid transit service ($F_P$-service). This design is justified by similar operations in Utah and Florida, as discussed in Section 1.

**A2** A many-to-one travel demand pattern is assumed here. The demand density function is exogenous and is expressed as a linear function, with the CBD demand density indicated by $Q_{CBD}$. Demand density decreases linearly, reaching zero at city boundary A.

**A3** It is assumed that in the peak periods, all commuters traveling along the corridor choose the lower travel cost mode. A logit model is used to capture the commuters' mode choice preference according to the travel cost comparison between auto and bus.

**A4** For the dynamic models, two stages are defined in Figure 4. The first stage, represented by $i = 0$, represents the stage in which fare-based transit service is implemented in the corridor. Conversely, stage 1, denoted by $i = 1$, refers to a $F_F$-service .

This dynamic adaptive operational approach is illustrated as follows. Initially, the transit agency may decide to run a $F_P$-service (or $F_F$-service) at time $t_0$. Later, it may decide to transition to a $F_F$-service (or $F_P$-service) at time $t_1$ when demand exceeds some threshold. This transition involves a switching cost, denoted by $D$. Later, depending on the evolution of the demand, the agency may decide to revert to a $F_P$-service from the $F_F$-service. This reversion comes with its own cost, labelled $K$, bringing the system back to the first stage, $i = 0$.

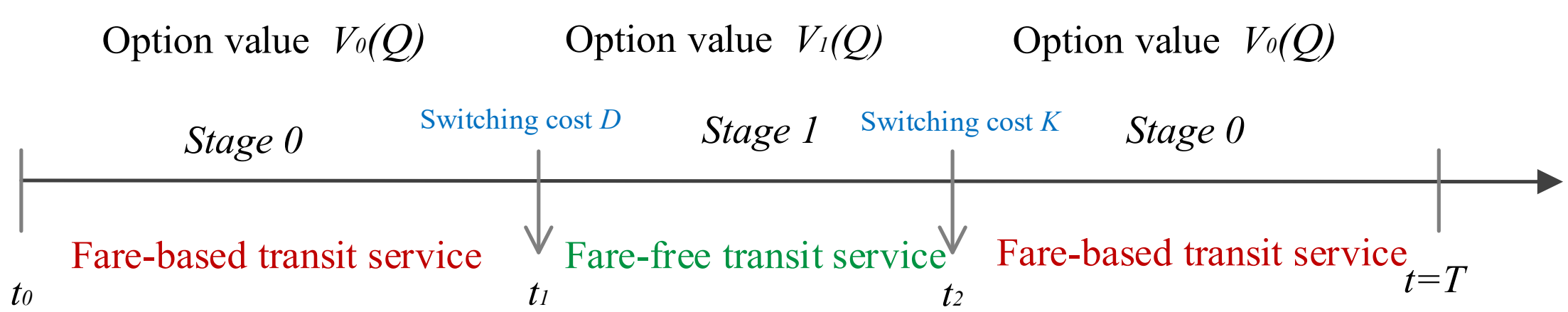

Figure 4: Fare-free transit service sequential implementation stages.

## 3.2 Static model setting

We start with the most basic policy and the model developed for analyzing it.

### 3.2.1 Travel demand and mode choice

Let $Q(x)$ be the potential density of travel demand (including auto and bus) at location $x$. $Q_{CBD}$ represents the demand density at CBD area, which is related via Eq. (1).

$$Q(x) = Q_{CBD} - \frac{Q_{CBD}}{A}x, \quad \forall x \in [0, A], \tag{1}$$

Referring to Figure 1, the transportation corridor can be covered by two different zones of transit service, namely the fare-free zone $[0, B]$ and fare-base zone$[B, A]$.

The user cost of traveling by bus from location $x$ to the CBD includes the in-vehicle travel time, waiting time, fare, and access time, which are expressed in Eq. (2).

$$C_b(x) = \begin{cases} \alpha_T \frac{x}{V_b} + \frac{1}{2F_i}\alpha_w + f_i + \alpha_S S, & \forall\, 0 \le x < B, \\ \alpha_T \frac{x}{V_b} + \frac{1}{2F_i}\alpha_w + f + \alpha_S S, & \forall\, B \le x \le A, \end{cases} \quad (2)$$

where $F_i$ is the transit frequency per hour at stage $i$. The average waiting time per trip is assumed to be half of the headway, which can be expressed as $\frac{1}{2F_i}$. $V_b$ is the average bus travel speed. The values of traveling time, waiting time and access time to and from stations are denoted as $\alpha_T$, $\alpha_w$, and $\alpha_S$, respectively. $S$ is the average access time to the nearest bus station. $f$ is the existing transit fare, which is given, and $f_i$ is the transit fare (where $f_0 = f, f_1 = \$0$) at different stages $i$.

The cost of traveling by auto from location $x$ to the CBD includes the in-vehicle travel time, fixed cost (e.g., parking charge in the CBD per work trip), and variable cost (e.g., auto gas price per unit of distance), which is expressed as:

$$C_a(x) = \alpha_T \frac{x}{V_a} + C_f^a + C_v^a x, \quad \forall 0 \le x \le A, \quad (3)$$

Given the transit service frequency and fare, the modal split can be formulated with a logit model as follows. The probability of selecting bus mode is:

$$Pr_b(x) = \begin{cases} \frac{exp(-\psi C_b(x))}{exp(-\psi C_a(x)) + exp(-\psi C_b(x))}, & \forall\, 0 \le x < B, \\ \frac{exp(-\psi C_b(x))}{exp(-\psi C_a(x)) + exp(-\psi C_b(x))}, & \forall\, B \le x \le A, \end{cases} \quad (4)$$

The probability of selecting bus in Eq. (4) can be rearranged as:

$$Pr_b(x) = \begin{cases} \frac{1}{exp(-\psi(\chi_0 + \chi_1)) + 1}, & \forall\, 0 \le x < B, \\ \frac{1}{exp(-\psi(\chi_0 + \chi_2)) + 1}, & \forall\, B \le x \le A, \end{cases} \quad (5)$$

where $\chi_0$, $\chi_1$ and $\chi_2$ can be determined by the difference between Eqs. (3) and (2), specifically by $\chi_0 = \alpha_T \left(\frac{1}{V_a} - \frac{1}{V_b}\right) + C_v^a$, $\chi_1 = C_f^a - \frac{1}{2F}\alpha_w - f_i - \alpha_S S$, and $\chi_2 = C_f^a - \frac{1}{2F}\alpha_w - f - \alpha_S S$.

In addition, the probability of selecting the auto mode is:

$$Pr_a(x) = \begin{cases} 1 - Pr_b(x), & \forall 0 \le x < B, \\ 1 - Pr_b(x), & \forall B \le x \le A, \end{cases} \quad (6)$$

The travel demand density by bus/auto at location $x$ can be expressed as the product of the total demand density and the probability of choosing the bus/auto. Let $Q_m(x)$ denote the travel demand density by bus/auto at location $x$, which can be expressed as:

$$Q_m(x) = Q(x) Pr_m(x), \quad (7)$$

Let $q_n^m(x)$ be the actual demand density that is elastic with generalized cost $C^m(x)$.

$$q_m(x) = Q_m(x)\big(1 - e_c C_m(x)\big), \quad \forall x \in [0, A], \quad (8)$$

where $e_c$ is the coefficient of generalized cost (i.e., travel impedance) in the demand function. When travel impedance reaches $\frac{1}{e_c}$, no trips are generated.

### 3.2.2 Fare-free zone determination when maximizing social welfare

We can now focus on the optimization of a bimodal corridor system when bus services are provided by a public agency. The objective of the transit agency is to maximize the total social welfare of the system by simultaneously determining the zone length $B$ and bus frequency $F_i$ when the demand density at the CBD is given.

The total social welfare per year at binary stage $i$ is denoted as $W_i$, which is the sum of producer surplus and user surplus $W_{ui}$. The producer surplus is also called profit, which is the difference between fare revenue $R_i$, bus operation cost $W_{oi}$, and fare collection cost $W_{fi}$. The total social welfare $W_i$ is expressed as:

$$W_i = U_i + (R_i - W_{oi} - W_{fi} - W_{d1}), \forall i = 0, 1, \tag{9}$$

where $W_0 = U_0 + (R_0 - W_{o0} - W_{f0})$ implies that the total social welfare at stage 0 includes user surplus, operating costs and costs associated with fare collection for the entire transportation corridor. $W_1 = U_1 + (R_1 - W_{o1} - W_{f1} - W_{d1})$ implies that the total social welfare at stage 1 includes user surplus, operating costs and costs associated with fare collection specifically for the region between $B$ and the city boundary $A$.

The user surplus $U_i$ can be computed by integrating the inverted function over demand, as in Chang and Schonfeld (1993).

$$U_i = \int_{C_m(F_i, B_i)}^{\frac{1}{e_c}} q_m(x) dC_m(x), \quad \forall x \in [0, A], \tag{10}$$

After inserting Eq. (2) into Eq. (3), the user surplus becomes a function of *x* as follows:

$$U_i = Q_m(x) \left( \frac{1}{2e_c} - C_m(x) + \frac{e_c}{2} C_m(x)^2 \right), \quad \forall x \in [0, A], \tag{11}$$

The bus operating cost $W_{oi}$ consists of fixed and variable operating costs, which is expressed as:

$$W_{oi} = g_f + g_v N_i, \; \forall i = 0, 1, \tag{12}$$

where $N_i$ is the fleet size on the transit corridor at stage $i$, and $g_f$ and $g_v$ are the fixed and the variable operation cost parameters, respectively. The fleet size $N_i$ can be further determined as the product of the round-trip time and the frequency $F_i$ in Eq. (13).

$$N_i = \frac{2A}{V_b} \times F_i, \forall i = 0, 1, \tag{13}$$

The fare collection cost at stage 0, $W_{f0}$, is related to the entire route length A, which can be expressed as:

$$W_{f0} = e_0 + e_1 A + e_2 \int_0^A q_b(x) dx, \tag{14}$$

The fare collection cost at stage 1, $W_{f1}$, is related to the length of the fare charging zone *A-B* after implementation of $F_F$-service, which can be expressed as:

$$W_{f1} = e_0 + e_1 (A - B) + e_2 \int_B^A q_b(x) dx, \tag{15}$$

where $e_0$ represents the fixed cost associated with the activation of the fare collection system, $e_1$ is the variable cost per unit of length of fare collection system, and $e_2$ is the cost per trip.

The cost related to fare collection is eliminated, reducing overhead. However, new administrative costs arise for managing the expanded fare-free system, such as route planning, maintaining service quality, and managing larger ridership. As the fare-free zone expands, more resources are needed to ensure proper operation and management (Regional Zero-Fare Study Final Report). This includes staff costs, security measures, and infrastructure investments (e.g., fare-free zone setting) (Kirschen et al., 2022). Therefore, the administrative costs include the fixed and variable components as follows:

$$W_{d1} = \iota_f + \iota_v B^{\theta}, \tag{16}$$

where $\iota_f$ and $\iota_v$ represents the fixed cost and variable cost associated with the administration costs of fare-free program, and $\theta$ is a parameter reflecting the nonlinear nature of the cost growth, related to the fare-free zone length.

The social welfare maximization problem at stage 0, $W_0$, is formulated as the sum of user surplus and producer surplus:

$$\max_{\{F_0\}} W_0 = \int_0^A \big(U_b(x) + U_a(x)\big)dx + \int_0^A q_b(x) f dx - \left(g_f + \frac{2\, g_v A F_0}{V_b}\right) - \Big(e_0 + e_1 A + e_2 \int_0^A q_b(x)dx\Big), \quad \forall x \in [0, A], \tag{17}$$

The scalar optimized transit frequency $F_0$ can be found by numerical search of Eq. (17), which will be shown in the following section.

When the transit service frequency $F_0$ is given, the social welfare can be rearranged as a function of demand density at the CBD in Eq. (18) when the parameters $d_0^Q$ and $d_0^C$ are determined first.

$$W_0(Q, F_0) = d_0^Q(F_0)Q + d_0^C(F_0), \tag{18}$$

where $d_0^Q$, and $d_0^C$ are function of $F_0$ that can be determined by Eqs. (19) and (20).

$$d_0^Q(F_0) = \int_0^A (1 - x/A)\left(Pr_a(x)\left(\frac{1}{2e_c} - C_a(x) + \frac{e_c}{2} C_a(x)^2\right) + Pr_b(x)\left(\frac{1}{2e_c} - C_b(x) + \frac{e_c}{2} C_b(x)^2 + \big(1 - e_c C_b(x)\big)(f - e_2)\right)\right) dx, \tag{19}$$

$$d_0^C(F_0) = -g_0 - \frac{2\, g_v A F_0}{V_b} - e_0 - e_1 B, \tag{20}$$

Therefore, once the spatial distribution of demand along the corridor and the various cost values are given, we can associate an optimized transit frequency $F_0$ at stage 0.

The social welfare maximization problem at $W_1$ stage 1 is formulated as Eq. (21):

$$\max_{\{B, F_1\}} W_1 = \int_0^B \big(U_b(x) + U_a(x)\big)dx + \int_B^A \big(U_b(x) + U_a(x)\big)dx + \int_B^A q_b(x) f dx - \left(g_f + \frac{2\, g_v A F_1}{V_b}\right) - W_{r1} - \Big(e_0 + e_1(A - B) + e_2 \int_B^A q_b(x)dx\Big), \tag{21}$$

The optimized transit frequency $F_1$ and the $F_F$-service coverage area $B$ can be found through bidimensional numerical search of Eq. (21), which will be shown in the following section.

When the transit frequency $F_1$, and the zone length $B$ are given, the social welfare function can be rearranged as Eq. (22), which is a function of the demand density at the CBD when the parameters $d_1^Q$, and $d_1^C$ are determined first.

$$W_1(Q, B, F_1) = d_1^Q(B, F_1)Q + d_1^C(B, F_1), \tag{22}$$

where $d_1^Q$ and $d_1^C$ are functions of $B$ and $F_1$ which can be determined with Eqs. (23) and (24).

$$d_1^Q(B, F_1) = \int_0^B (1 - x/A)\left(Pr_a(x)\left(\frac{1}{2e_c} - C_a(x) + \frac{e_c}{2} C_a(x)^2\right) + Pr_b(x)\left(\frac{1}{2e_c} - C_b(x) + \frac{e_c}{2} C_b(x)^2\right)\right) dx + \int_B^A (1 - x/A)\left(Pr_a(x)\left(\frac{1}{2e_c} - C_a(x) + \frac{e_c}{2} C_a(x)^2\right) + Pr_b^H(x)\left(\left(\frac{1}{2e_c} - C_b(x) + \frac{e_c}{2} C_b(x)^2\right) + \big(1 - e_c C_b(x)\big)(f - e_2)\right)\right) dx, \tag{23}$$

$$d_1^C(B, F_1) = -g_f - \frac{2\, g_v A F_1}{V_b} - \iota_f - \iota_v B^\theta - e_0 - e_1(A - B). \tag{24}$$

Therefore, whenever the demand densities at the CBD $Q_{CBD}$ is given, we can associate the optimized transit frequency $F_1$ and the zone length $B$ at stage 1.

The instantaneous social welfare improvement can be obtained in Eq. (25) from the difference between social welfare at stage 1 from Eq. (22) and stage 0 from Eq. (18):

$$\Omega(Q) = W_1(Q) - W_0(Q) = \big(d_1^Q - d_0^Q\big)Q + \big(d_1^C - d_0^C\big). \tag{25}$$

### 3.2.3 Social equity

The previous research on $F_F$-service has frequently overlooked an intriguing observation: the implementation of $F_F$-service can result in significant variations in user surplus among commuters along the transportation corridor (Wang et al., 2021). This implies that individuals from diverse demand groups along the corridor may experience dissimilar advantages or disadvantages as a result of the $F_F$-service. Such a situation indicates a lack of social equity (Farber et al., 2014; Zhou et al., 2019). Therefore, it is crucial for decision-makers to consider the equity aspect when implementing $F_F$-service.

According to Assumption 2, the spatial demand is illustrated in Figure 3, where the demand is divided into n groups from the CBD to city boundary as ($P_1$, $P_1, P_{...}$, $P_n$). It is assumed that the user surplus per person within each group $P_\zeta$ is the same, denoted as $\overline{U_{ab}}^{\zeta}$. The Gini index is used to define the level of inequality in fare-free transit, as outlined in the subsequent text. The Gini index is a quantitative measure of inequality, with values ranging from 0 to 1, and has been applied in several transportation fields (Li & DaCosta, 2013; Chen et al., 2019; Ben-Elia & Benenson, 2019; Zhang et al., 2020; Gao& Li, 2024). In the context of user surplus distribution, the Gini coefficient measures the fairness of surplus allocation across different groups along the transportation corridor.

**Definition 1:** Let $E_i(B)$ denote the inequality ratio of user surplus for demand groups distributed along transportation corridor with Gini index:

$$E(B) = \frac{\sum_{\zeta=0}^{n}\sum_{\varrho=0}^{n}|x_\zeta - x_\varrho|}{2*n\sum_{\zeta=0}^{n}x_\zeta} = \frac{1}{2n^2\sum_{\zeta=0}^{n}\overline{U_{ab}}^{\zeta}}\sum_{\zeta\in n}\sum_{\varrho\in n}\left|\overline{U_{ab}}^{\zeta} - \overline{U_{ab}}^{\varrho}\right|, \tag{26}$$

where $n$ is the number of demand groups spatially distributed along the corridor extending from the CBD to the city boundary; $\overline{U_{ab}}^{\zeta}$ represents the average user surplus of demand groups $\zeta$ including bus and auto mode can be determined by Eq. (29).

$$\overline{U_{ab}}^{\zeta} = \frac{\int_{n_\zeta}^{n_{\zeta+1}}(U_b(x)+U_a(x))dx}{\int_{n_\zeta}^{n_{\zeta+1}}Q(x)dx}, \tag{27}$$

$E(B) \to 0$ indicates that average user surplus of different groups are even distributed across different groups, meaning that all demand groups benefit equally from implementing $F_F$-service. On the other hand, $E(B) \to 1$ signifies that the average user surplus of different demand groups are distributed most skewedly, with certain groups benefiting disproportionately, leading to extreme inequality for the user surplus along the traffic when implementing $F_F$-service. However, the values of 0 and 1 represent theoretical extremes that are unlikely to occur in real-world scenarios. Consequently, evaluation of Gini index can be defined as the degree of equity remains within acceptable bounds. By addressing this concern, decision-makers can strive to implement $F_F$-service schemes (coverage length $B$) that minimize any resulting unfairness.

### 3.2.4 Benefit index

In this section, a social benefit index will be defined to consider both social welfare and social equity for the length of $F_F$-service. After determining the optimized length of $F_F$-service from Section 3.2.2, a social welfare index of different fare-free zone lengths can be developed by normalizing the socially welfare from $F_P$-service to $F_F$-service. The social welfare index is expressed as:

$$SW(B) = \frac{W_1(B)-W_0}{W_1(B^*)-W_0}, \tag{28}$$

A novel benefit measure combines social welfare and social equity, which balances two competing objectives, namely social welfare and social equity, in relation to the fare-free coverage zone length. It is expressed as:

$$Max: G(B) = (1-\mu)\cdot SW(B) + \mu\cdot\big(1-E(B)\big), \tag{29}$$

Subject to:

$$\beta_{min} \le E(B) \le \beta_{max}. \tag{30}$$

where $\mu$ is a weighting factor (between 0 and 1) that determines the relative importance placed on social equity and social welfare. When $\mu$ approaches 1, greater emphasis is placed on equity; when $\mu$ approaches 0, more priority is given to social welfare. In social equity a higher Gini Index indicates greater inequality, while our goal in the benefit index is to maximize the benefit. To align the objective of maximizing benefit while considering both social welfare and social equity, we define social equity as 1 - Gini Index in benefit index. This transformation allows us to optimize social welfare and social equity in the same direction. Eq. (30) ensures the equity ratio under different scenarios of fare-free transit schemes for different demand groups does not exceed the equity level between $\beta_{min} = 0$ and $\beta_{max} = 1$. In other words, an equitable $F_F$-service requires that benefits or losses of demand groups should be restricted within an acceptable interval.

## 3.3 Dynamic model setting

For uncertain demand, a dynamic model is proposed to optimize the timing for switching the fare policy. Within this framework, Geometric Brownian motion (GBM) is employed to capture the random fluctuations in demand, as commonly used in financial modeling due to its capacity to simulate the stochastic behavior of financial assets (Giversen & Bendkia, 2011). It is further assumed that the demand density at the CBD $Q_t$ stochastically fluctuates over time as a non-stationary stochastic process. A GBM model is used to capture the overall demand density at the CBD $Q_t$ at time *t*, as shown in Eq. (31):

$$\frac{dQ_t}{Q_t} = \eta dt + \sigma dw(t), \tag{31}$$

where $\eta$ is the CBD total demand growth rate, $dt$ is an infinitesimal time increment, and $\sigma$ is the volatility rate which models the extent to which demand randomly changes. $dw(t)$ is an increment of a standard Wiener process. By definition, $w(t) = \varepsilon_t\sqrt{t}$ , where $\varepsilon_t$ is a standard normal random variable, whose mean is 0 and standard deviation is 1.

According to Merton (1976), the solution of Eq. (32) is:

$$Q_t = Q_{CBD}e^{(\eta-\frac{\sigma^2}{2})t+\sigma w(t)}. \tag{32}$$

where $Q_{CBD}$ is the present demand density at the CBD (when time *t = 0*).

Given Eq. (32), the expectation of $Q_t$ for any given value of $t$ can be derived as follows (Ross 1996):

$$\mathbf{E}[Q_t] = Q_{CBD}exp\,(\eta t). \tag{33}$$

We consider a stochastic process combined with optimal stopping in order to discover the form of the optimal strategy. We denote as $\mathscr{H}$ the family of adapted, finite variation, càdlàg process $H$ with value of $\{0, 1\}$ and by $\mathscr{L}$ the set of all stopping times. According to the market entry-exit real option (Zervos, 2003), a stochastic system is formulated that can operate in two states, an "active" one ($F_F$-service) and an "inactive" one ($F_P$-service). The system's operation status can be changed at a sequence of stopping times. These transitions times constitute a decision policy that we model by a process $H \in \mathscr{H}$. Specifically, given any time *t*, $H_t = 1$=1 if the system is "active" ($F_F$-service) at time *t*, whereas $H_t = 0$=0 if the system is "inactive" ($F_P$-service) at time *t*. The problem is to optimize the implementation strategy that maximizes the expected total payoff over the planning horizon *T*. The expected payoff is equal to realized social welfare improvement between two modes minus discounted switching cost occurred. Decisions on service mode are assumed to be activated instantaneously, i.e., with negligible delay, so with respect to the time horizon:

$$\mathcal{I}_{H,Q}(\mathcal{H},T) = \sup_{(\mathcal{H},T)} \boldsymbol{E}\left[\int_0^T [W_0(Q_s)H_s + W_1(Q_s)(1-H_s)]ds - \sum_{0\le s\le T} O_s[D(\Delta H_s)^+ + K(\Delta H_s)^-]\right]. \quad (34)$$

where the term $W_0(Q_s)H_s + W_1(Q_s)(1-H_s)$ represents the realized social welfare of two modes. The second term represents the expected occurred switching cost. $\Delta H_t = H_{t+} - H_{t-}$, $(\Delta H_s)^{\pm} = max\{\pm\Delta H_t, 0\}$, and the discounting process $O_s$ is given by *exp(-kt)*. *D, K* are the switching costs associated, respectively, with transition from its "active" mode to its "inactive" mode, and vice versa, at time t.

Specifically, we consider operating $F_P$-service ($V_0(Q)$) and activating $F_F$-service by ($V_1(Q)$) as the option value, whose value varies with demand density $Q$. The asset equilibrium condition, per Ito's lemma (Dixit, 1994), is equivalent to the second order differential equation for the transit service with fare described in Eq. (35). Similar computations can be used to determine the option value for activation of $F_F$-service, which requires considering the additional immediate social welfare improvement of implementing $F_F$-service. The ordinary differential equation can be solved to obtain the value of $V_1(Q)$ that satisfies the asset equilibrium condition for activating of $F_F$-service, as defined in Eq. (36).

$$\frac{1}{2}\sigma^2 Q^2 V_0''(Q) + \eta Q V_0'(Q) - kV_0(Q) = 0, \quad (35)$$

$$\frac{1}{2}\sigma^2 Q^2 V_1''(Q) + \eta Q V_1'(Q) - kV_1(Q) + \Omega(Q) = 0. \quad (36)$$

where $V_0(Q)$ is option value of operating $F_P$-service, $V_1(Q)$ is option value of activation of $F_F$-service, and $\Omega(Q)$ is the instantaneous social welfare improvement given in Eq. (25).

By accounting for ideal transitions between $F_P$-service and $F_F$-service, the best activation and cancellation timing policy is identified. The threshold for switching from one service to another is not unique because switching costs exist. Instead, switching costs yield two switching thresholds, which we treat as the upper bound demand threshold $\overline{Q}$ (from no fares to fares) and the lower bound demand threshold $\underline{Q}$ (from fares to no fares), respectively, for "entry and exit" (switching from $F_P$-service to $F_F$-service) with fare is treated as a market entry and vice versa for market exit) for the service switching.

The following "value matching" interactions between the option values are created by these thresholds, as defined in Eqs. (37) and Eq. (38):

$$V_0(\overline{Q}) = V_1(\overline{Q}) - D, \quad (37)$$

$$V_1\left(\underline{Q}\right) = V_0\left(\underline{Q}\right) - K. \quad (38)$$

where $D$ is the switching cost assumed for switching from no fares to fares, $K$ is the switching cost assumed for switching from $F_P$-service to $F_F$-service.

The smooth-pasting conditions can be expressed by Eqs. (39) - (40):

$$V_0'(\overline{Q}) = V_1'(\overline{Q}), \quad (39)$$

$$V_0'\left(\underline{Q}\right) = V_1'\left(\underline{Q}\right). \quad (40)$$

When demand density $Q$ is assumed to grow as a GBM process, as in Eq. (32), we have a general solution of $V_0(Q)$ using infinite series (Dixit, 1989), as shown later in Eq. (41).

$$V_0(Q) = X_0 Q^{\gamma_0} + Y_0 Q^{\gamma_1}. \quad (41)$$

Similarly, for $V_1(Q)$ the solution for Eq. (38) is obtained in Eq. (42):

$$V_1(Q) = X_1 Q^{\gamma_0} + Y_1 Q^{\gamma_1} + \Omega(Q). \quad (42)$$

where $X_0$, $Y_0$, $X_1$ and $Y_1$ are the constants to be determined, and $\Omega(Q)$ is given by Eq. (25).

$\gamma_0$ and $\gamma_1$ are the negative and positive root of the quadratic equation in Eqs. (44) - (45):

$$\frac{1}{2}\sigma^2\gamma(\gamma - 1) + \eta\gamma - k = 0. \quad (43)$$

According to the quadratic formula, they can be expressed by following (Dixit, 1989):

$$\gamma_1 = \frac{1}{2} - \frac{\eta}{\sigma^2} + \sqrt{\left(\frac{1}{2} - \frac{\eta}{\sigma^2}\right)^2 + \frac{2k}{\sigma^2}} > 1, \quad (44)$$

$$\gamma_0 = \frac{1}{2} - \frac{\eta}{\sigma^2} - \sqrt{\left(\frac{1}{2} - \frac{\eta}{\sigma^2}\right)^2 + \frac{2k}{\sigma^2}} < 0. \quad (45)$$

**Proposition 1**. When the switching costs between no fares and fares ($D$ and $K$) are positive, the upper and lower thresholds can be obtained as:

$$\left(\frac{\overline{Q}}{\underline{Q}}\right)^{\gamma_1} = \frac{-(\nu_1 - D)\gamma_0 + \nu_0(\gamma_0 - 1)\overline{Q}}{-(\nu_1 - K)\gamma_0 + \nu_0(\gamma_0 - 1)\underline{Q}} \quad (46)$$

$$\left(\frac{\overline{Q}}{\underline{Q}}\right)^{\gamma_0} = \frac{(\nu_1 - D)\gamma_1 + \nu_0(\gamma_1 - 1)\overline{Q}}{(\nu_1 - K)\gamma_1 + \nu_0(\gamma_1 - 1)\underline{Q}}. \quad (47)$$

where $\nu_0 = \frac{\left(d_1^Q - d_0^Q\right)}{k - \eta}$, and $\nu_1 = \frac{\left(d_1^C - d_0^C\right)}{k}$. From the system of non-linear equations (37) - (40), we can solve for the values of $\overline{Q}$ and $\underline{Q}$ numerically. The proof is given in Appendix A.

**Proposition 2.** When the switching cost between $F_F$-service and $F_P$-service become zero, the upper and lower demand thresholds converge to a single threshold $Q^*$, which can be expressed as follows:

$$Q^* = \frac{(\gamma_1 - \gamma_0)\gamma_0\gamma_1\nu_1}{\left((1-\gamma_0)\gamma_1(\gamma_1 - 1) - (1-\gamma_1)\gamma_0(\gamma_0 - 1)\right)\nu_0} = \frac{\gamma_0\gamma_1\nu_1}{(1-\gamma_0)(\gamma_1 - 1)\nu_0}. \quad (48)$$

The proof is given in Appendix B.

# 4. Experiments and Discussion

A series of increasingly complex scenarios are studied, beginning with an exploration of optimal design for $F_P$-service within the static model. By considering the given demand density at a particular time $t$, the optimized configuration can be represented in terms of zone length, transit frequency, and mode split. The crux of this optimal design resides in striking a balance between user surplus and producer surplus. The Gini-based equity index for demand groups along the corridor is examined under varying transit fares. Then, the Lorenz curve is used to identify the fare policy that offers the highest level of equity (Wang et al., 2021). By combining social welfare and social equity, a benefit index is analyzed, weighing social equity and social welfare differently.

Moving forward, the analysis is further extended to a real-world transit line with our proposed dynamic model, incorporating demand uncertainty through a GBM process. Within the realm of demand uncertainty, the potential implementation of $F_P$-service can be treated as an option-based optimal policy that simply depends on the current state, which is analytically calculable. The evaluation then focuses on the policies for implementation (timing of activation and de-activation) of $F_P$-service, and derive the optimized zone length, transit frequency, and associated social equity ratio over time.

Lastly, verification analyses are conducted to assess the computational performance of various policies, such as the deterministic dynamic policy, stochastic policy without switching costs, stochastic policy with switching costs, retaining the existing $F_P$-service, and adhering to $F_F$-service.

## 4.1 Data input

Unless otherwise noted, the input parameters and baseline values in Table 2 are mostly based on Li et al. (2015); Guo et al. (2021) and Tirachini and Hensher (2011). Other values of some key parameters are used in sensitivity analyses in Section 4.4.

Table 2: Baseline values

| Notation | Definition | Baseline value | Unit |
|---|---|---|---|
| $n$ | The number of demand groups spatially distributed along the corridor extended from the CBD to the city boundary | 50 | - |
| $Q_{CBD}$ | Demand density at CBD point | 1500 | pax/mile/day |
| $\psi$ | The parameter in the Logit assignment model | 0.5 | - |
| $A$ | City boundary | 50 | mile |
| $e_c$ | Coefficient of generalized cost of auto users | $\frac{1}{70}$ | - |
| $\alpha_T$ | Value of in-vehicle travel time | 20 | \$/hr |
| $\alpha_w$ | Value of waiting time | 20 | \$/hr |
| $C_f^a$ | Fixed cost of auto user | 3 | \$/trip |
| $C_v^a$ | Variable cost for auto user | 0.5 | \$/mile |
| $f_i$ | Transit fare at stage $i$ | 5 (or 0) | \$/trip |
| $S$ | Average access time to transit stops | 6 | mins |
| $\alpha_S$ | Value of access time | 20 | \$/hr |
| $g_f$ | Fixed cost for transit operation | 10,000 | \$/day |
| $g_v$ | Variable cost for transit operation | 500 | \$/vehicle/day |
| $e_0$ | Fixed cost associated with the establishment of the fare collection system | 5,000 | \$/day |
| $e_1$ | Variable cost per unit of length of fare collection system | 500 | \$/mile |
| $e_2$ | Fare collection cost per trip | 0.1 | \$/trip |
| $\iota_f$ | Fixed administration cost of fare-free program | 10000 | \$/day |
| $\iota_v$ | Variable administration cost of administration cost for fare-free program | 10 | \$/mile |
| $\theta$ | Parameters to reflect the nonlinear nature of the cost growth | 2 | - |
| $\eta$ | Monthly growth rate of CBD demand | 1.16 | % |
| $\sigma$ | Monthly volatility rate CBD demand | 13.47 | % |
| $k$ | Monthly discount rate | 0.02 | - |
| $D$ | Cost of switching from $F_F$-service to $F_P$-service | 5,000 | \$ |
| $K$ | Cost of switching from $F_P$-service to $F_F$-service | 5,000 | \$ |

## 4.2 Static models analysis

### 4.2.1 Optimized fare-free transit service design

Figure 5 demonstrates how the total social welfare of the $F_F$-service changes with respect to the fare-free zone length $B$ and transit frequency $F$. The optimized fare-free zone length and transit frequency for the $F_F$-service can be determined by searching for the total social welfare. Social welfare comprises both auto user and bus user surplus and various cost, which are expressed in Eq. (21). The contour lines form concentric shapes, with a peak point located at 32 miles from the CBD and a transit frequency of 13. The circle marks the most beneficial solution, namely $B^* = 32$ and $F^* = 13$ that yields the highest value. When the fare-free zone is fixed (e.g., at 32), increasing the transit frequency from 2 to 13 initially boosts social welfare. However, beyond this point, further increases in frequency lead to a decline in social welfare. This pattern arises because increasing transit frequency does offer benefits, such as reduced waiting times for passengers, but the marginal benefits diminish as the rising operating costs associated with higher frequency outweigh these advantages. Moreover, when examining the variation in total social welfare, changes in the fare-free zone length have a more gradual impact compared to adjustments in transit frequency. This means that social welfare is less sensitive to changes in the fare-free zone length than to changes in transit frequency.

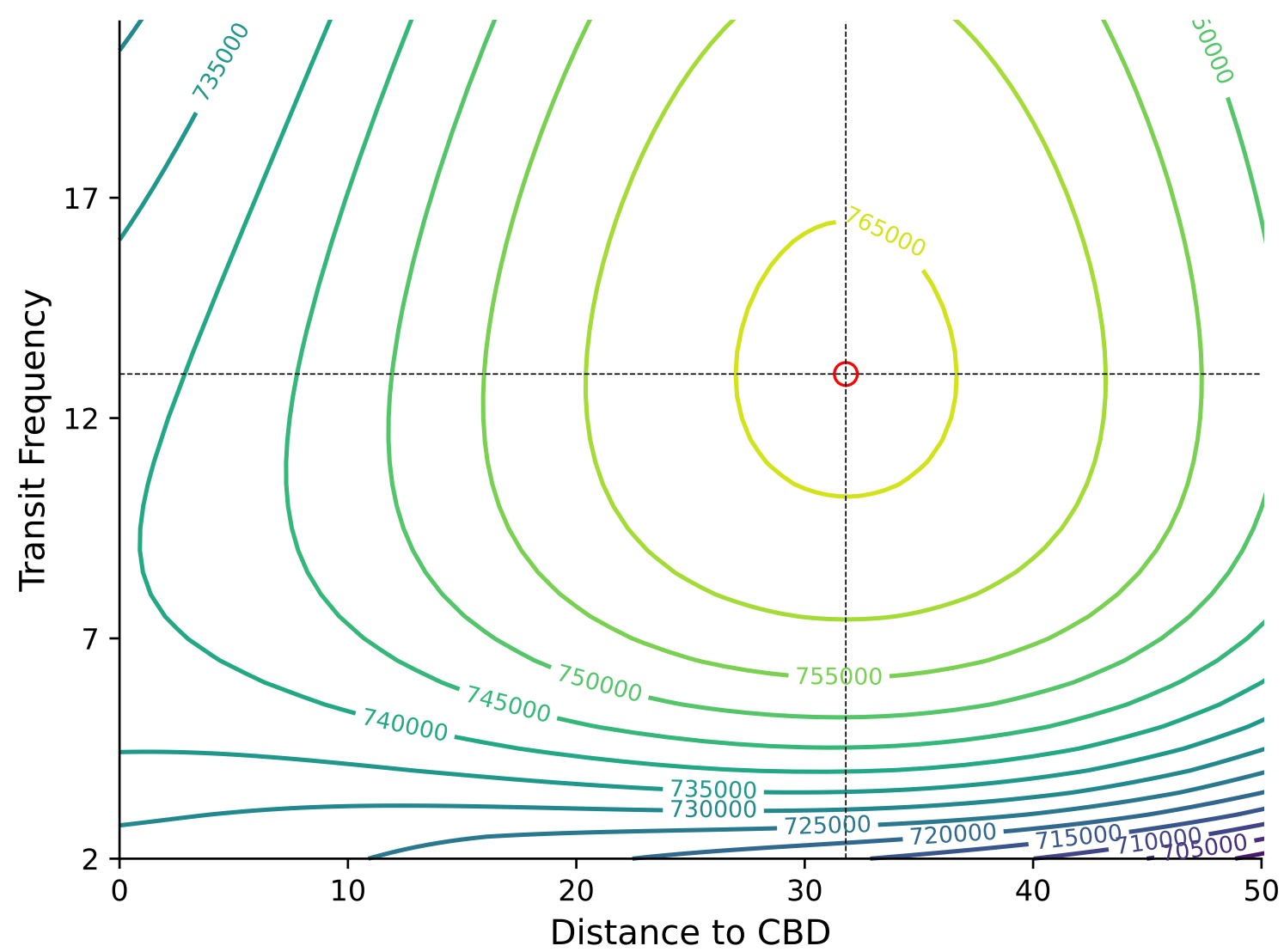

Figure 5: The optimized fare-free transit service zone length at different transit frequencies.

The population of bus and auto users can be observed in Figures 6 (a) and (b), respectively. It is notable that the population of both auto and bus users within the zone increases as the fare-free zone lengthens. Conversely, if the zone extends beyond a certain point, the population of auto and bus users outside the zone area decreases. Figure 7 presents the auto and bus mode shares when the zone length is preset at 32 miles from the CBD. The mode share curve appears disjointed due to the switch from zero to positive transit fares at the border of the fare-free zone. Bus share is concentrated primarily along transit corridors. The transit fare variation between the zone and outside the zone leads to a significant 40% shift in the mode share. Expanding the length of the fare-free zone stimulates more demand for bus trips within the zone, leading to a higher bus share. Therefore, outside of the fare-free zone, the bus share declines while the auto share increases.

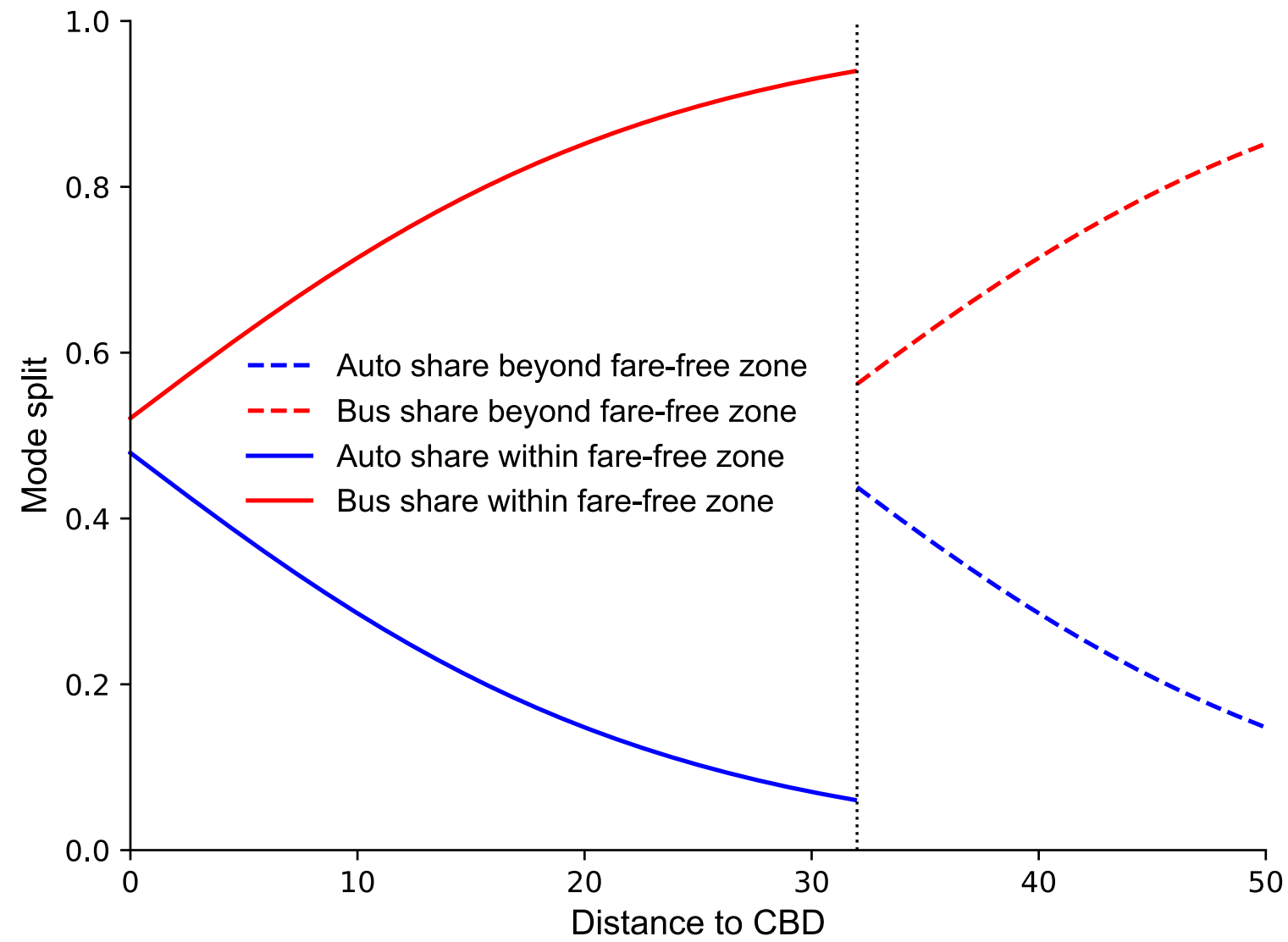

Figure 7: Frequency-optimized mode splits vs. the distance to the CBD when the length of fare-free zone is fixed at 32 miles

### 4.2.2 Social equity differences

Figure 8 illustrates the user surplus within and beyond fare-free zone based on the provided parameters. The user surplus is represented by four lines, depicting the bus and auto user surplus within and beyond the fare-free zone. The bus user surplus is indicated by dashed lines, while the auto user surplus is represented by solid lines. At 32 miles from the CBD, the bus user lines are disconnected because bus users are required to pay the fare outside the zone and vice versa. For auto users, the user surplus is also discontinuous within and outside the region due to the variation in auto share between these regions. Outside the region, bus user surplus declines while auto user surplus increases because the demand from bus users decrease. Outside of the region, the user surplus declines more for bus than for auto.

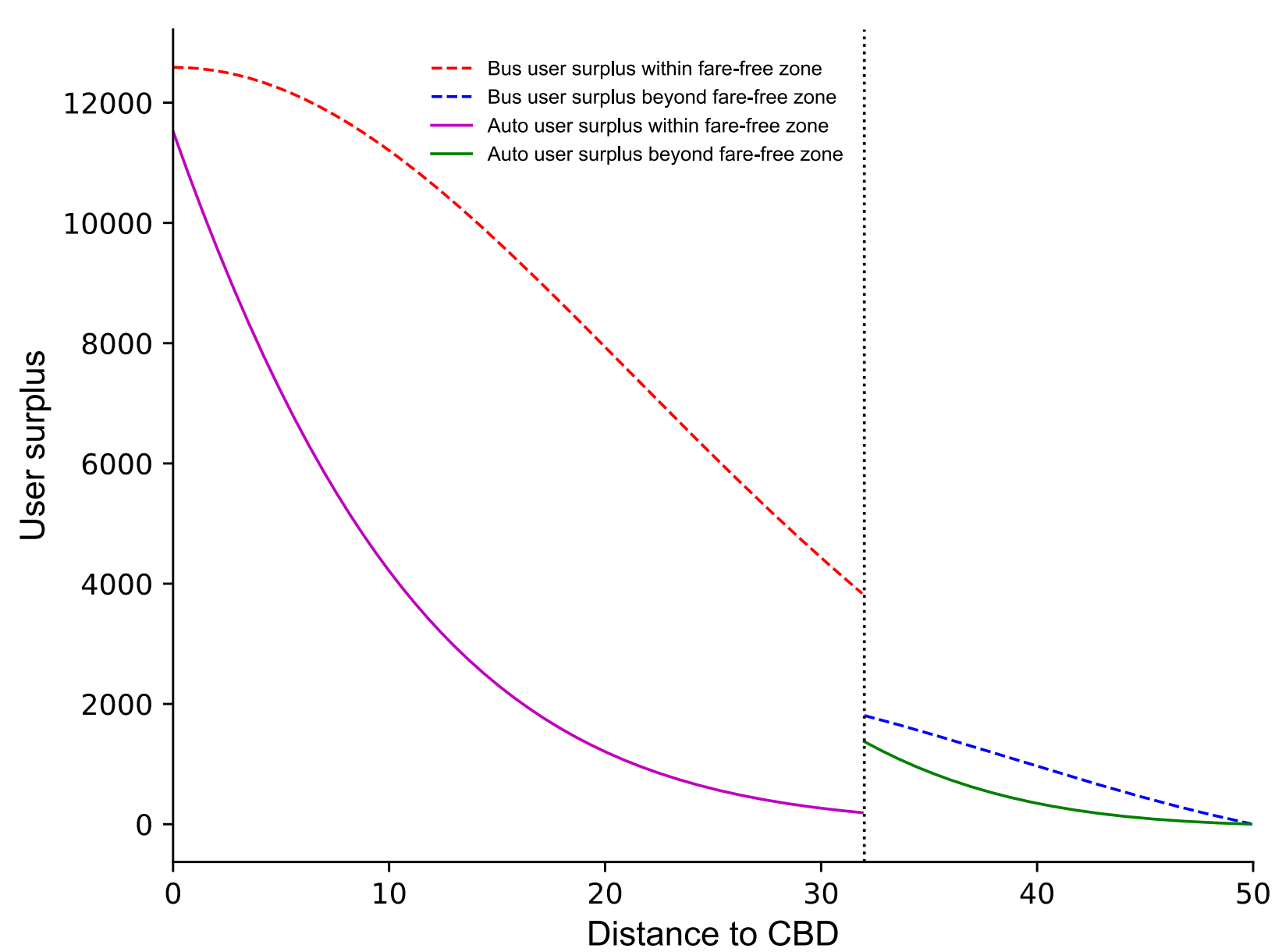

Figure 8: Frequency-optimized user costs vs. distance to the CBD by mode when $B$ = 32 mi.

When the user surplus values are provided as base values, the social inequality ratio (Gini index) can be computed with Eq. (26). This ratio represents the equity of transit fare policy for demand groups distributed in the corridor. Figure 9 illustrates how the Gini index changes in response to variations in transit fare and the length of the fare-free zone. When the fare is above \$0, the Gini index initially increases, reflecting a rise in inequality, but then gradually decreases, reaching its lowest level when the $F_F$-service coverage extends to the city boundary. As the transit fare decreases from \$8 to \$0, the Gini index decreases, which means that greater equity between demand groups along the corridor. This is because lower or eliminated fares make transit services more accessible to all groups, reducing the financial barriers that might otherwise limit transit usage. Consequently, the distribution of transit benefits becomes more balanced, reducing disparities between demand groups and promoting equity. Therefore, the highest level of equity for different demand groups (at the lowest Gini index) is reached when implementing $F_F$-service, where the fare is equal to \$0.

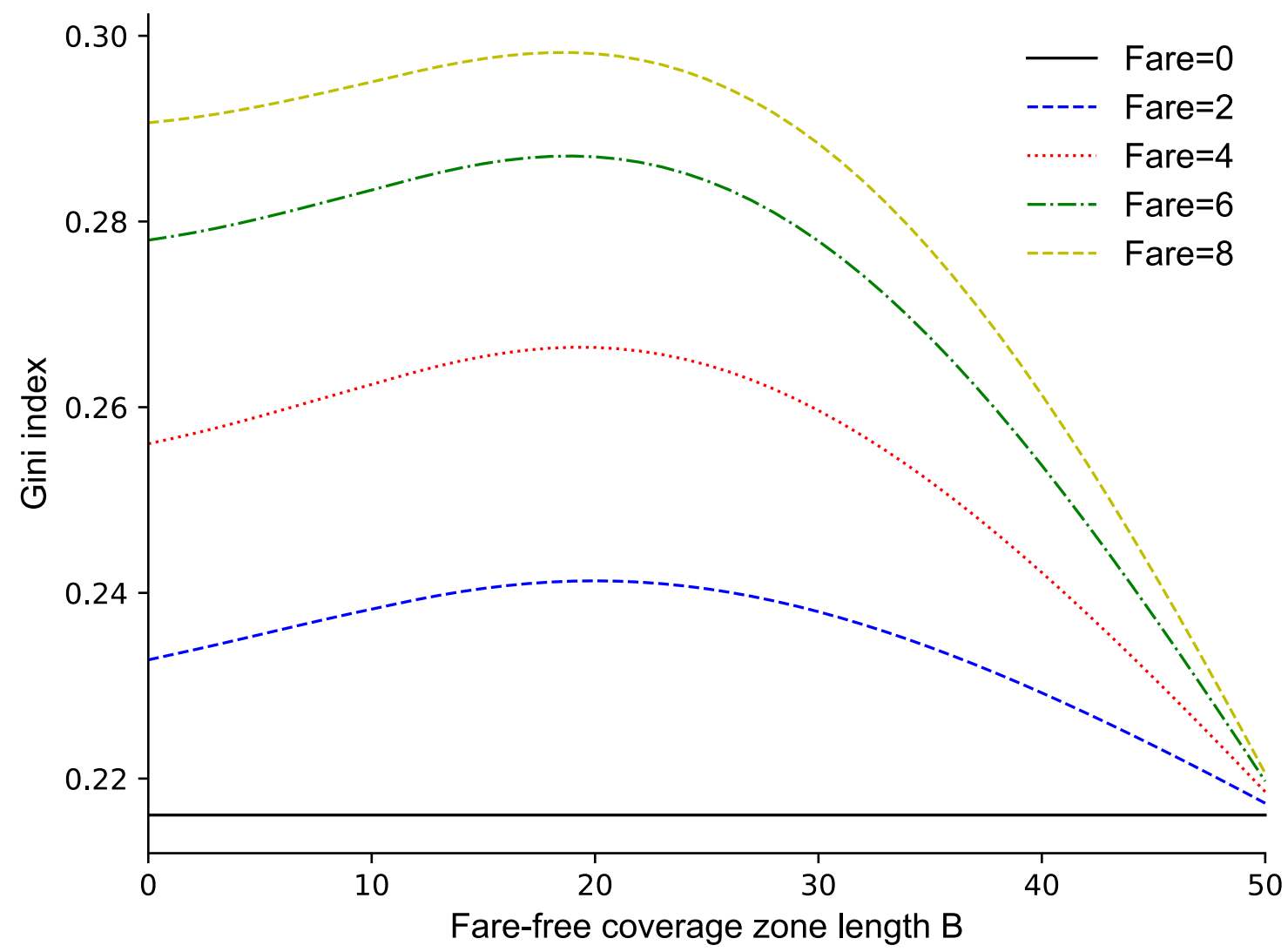

Figure 9: Inequality ratio (Gini index) of different demand groups vs. zone length at various fares.

We then evaluate equity using Lorenz curves and Gini indexes for different demand groups. The Lorenz curve is a standard economic tool used to measure equity based on the relative distribution within a population (Wang et al., 2021). To generate the Lorenz curve, we first rank the average user surplus of each demand group in ascending order and then plot the cumulative proportion of the travel demand on the horizontal axis and the cumulative proportion of average user surplus on the vertical axis. If the average user surplus across all demand groups is equal, the result is a black dotted line (the Line of Equity) in Figure 10, illustrating a linear relation between the average user surplus and demand distribution. However, in practice, the Lorenz curve deviates from this line of equity. The closer the curve is to the Line of Equity, the more equal the distribution (Delbosc & Currie, 2011).

In Figure 10, three Lorenz curves are generated to compare the equity for different demand groups of different fare policies among $F_P$-service, transit with optimized fare-free zone for $F_F$-service. The Lorenz curve of $F_F$-service in the whole corridor (green dash-dot line) is closer to the Line of Equity than the curve of $F_P$-service in whole corridor (blue dashed line) and transit service with optimized fare-free zone (orange solid line), indicating an improvement in equity when implementing the $F_F$-service in the whole corridor. This finding is also reflected in the decrease in Gini Index from 0.2623 and 0.2607 to 0.2113, since a smaller Gini Index means better equity.

Additionally, transit service with an optimized fare-free zone also improves the equity compared with fare-based transit in whole corridor, though the Gini index increases very slightly from 0.2603 to 0.2607. Therefore, the $F_F$-service in whole corridor is the most beneficial choice from the perspective of equity.

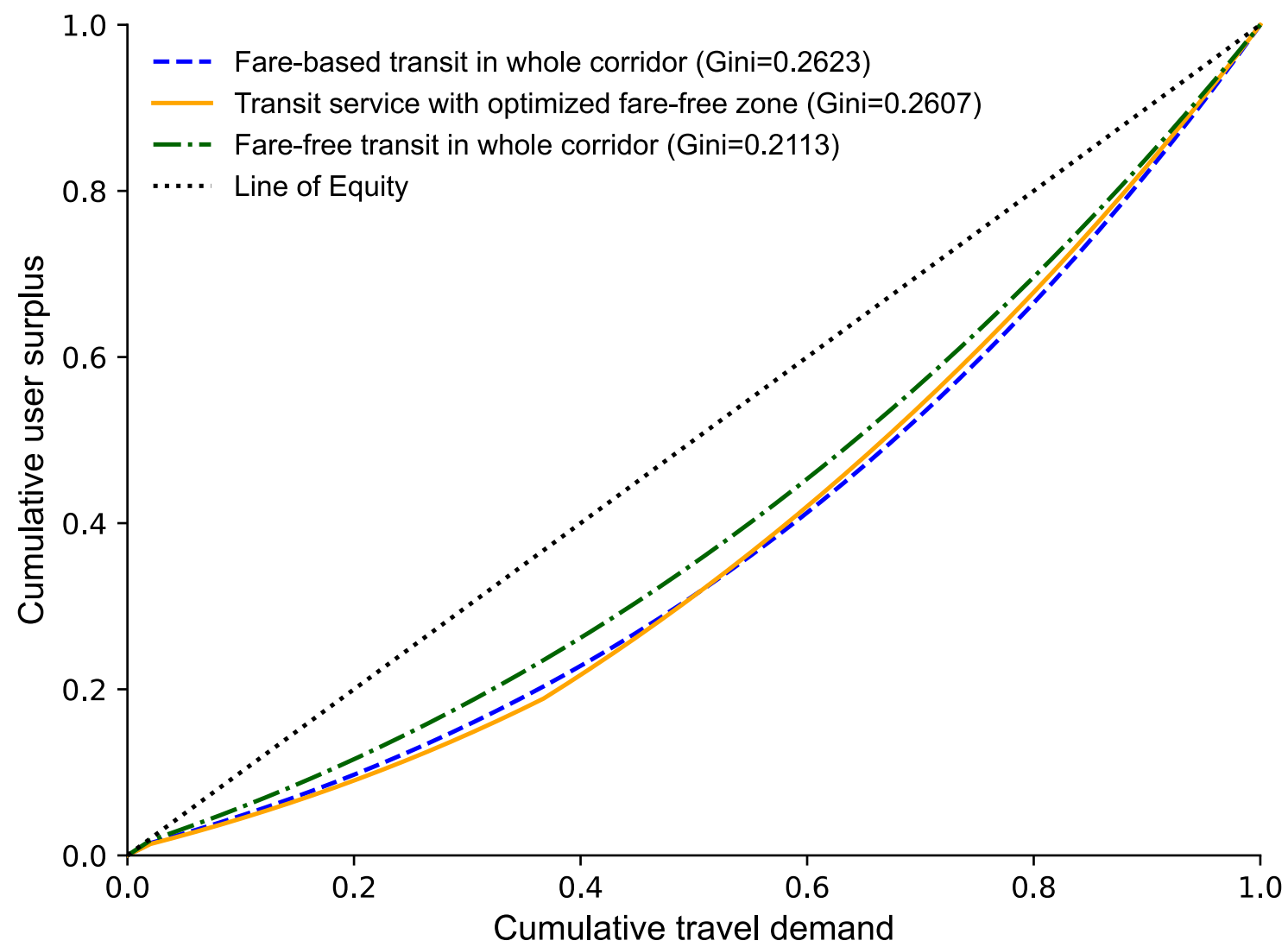

Figure 10: The level of equity across different fare policy.

### 4.2.3 Benefit measure analysis

Figure 11 illustrates the changes in normalized social welfare (dashed line) and social equity (dotted line) relative to the length of the fare-free zone. From this figure, the most beneficial fare-free zone length can be identified from two perspectives: social welfare and social equity. When considering only social welfare, the most favorable fare-free zone length is 32 (marked by the circle), as detailed in the previous section. On the other hand, when focusing solely on social equity, the highest level of equity occurs when the fare-free zone spans the entire 50 mile corridor. This is further compared with the Gini index in Figure 10. By combining social welfare and social equity, a benefit measure is introduced in section 3.2.4. When the parameter $\mu = 0.5$, the benefit measure is represented by a solid line in Figure 11.

Figure 12 shows how the benefit measure and the optimized length of the fare-free zone vary for different weight values of $\mu$. As $\mu$ increases, meaning that more weight is placed on social equity in the benefit measure calculation, the optimized fare-free zone length also increases, as a longer zone enhances social equity. However, since the length of the fare-free zone has more influence on social welfare than on social equity, smaller values of $\mu$ (e.g. $\mu = 0.1$, $\mu = 0.3$ and $\mu = 0.5$) have less effect on the optimized length of the fare-free zone.

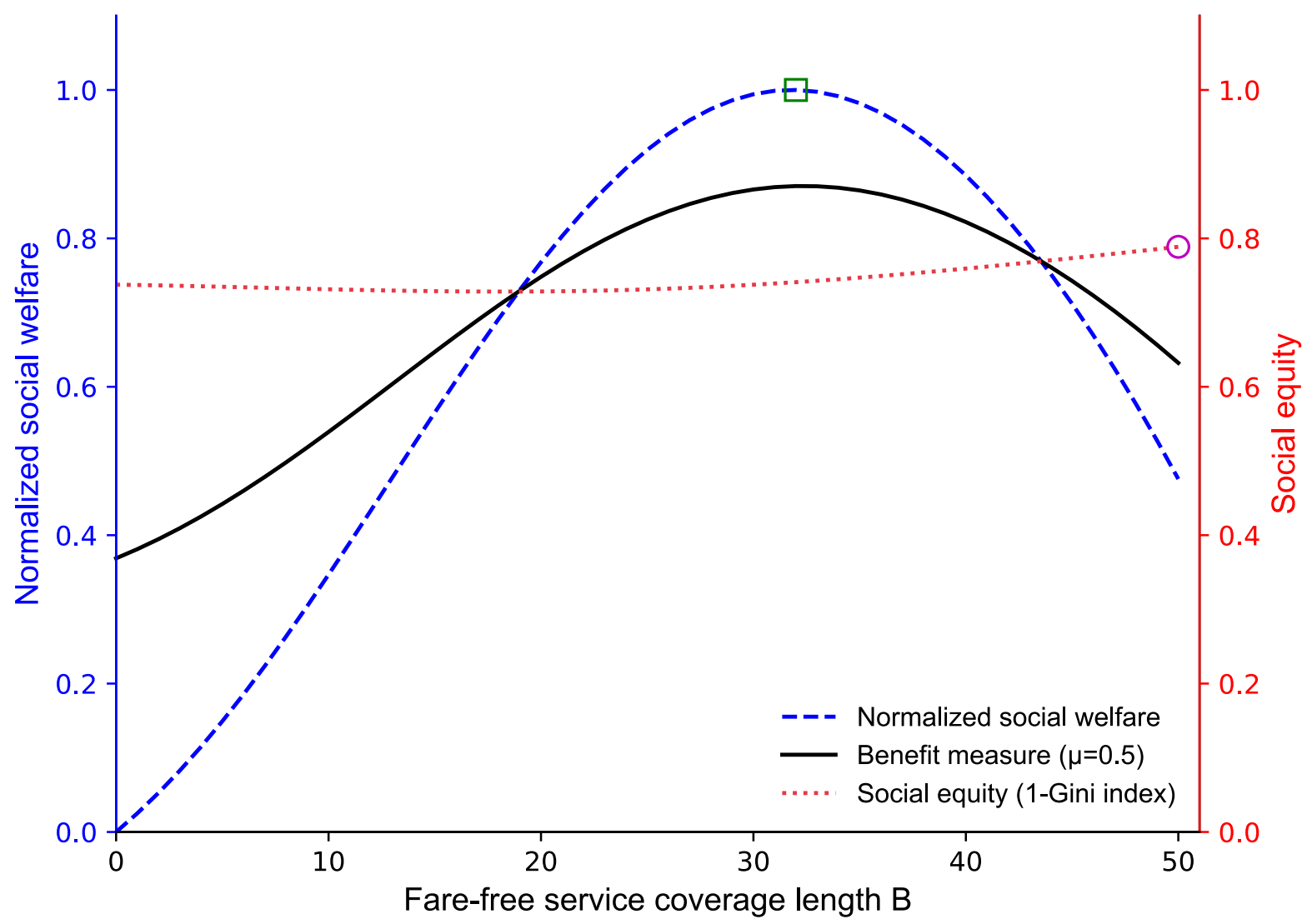

Figure 11: Optimized zone length and equity-constrained zone length.

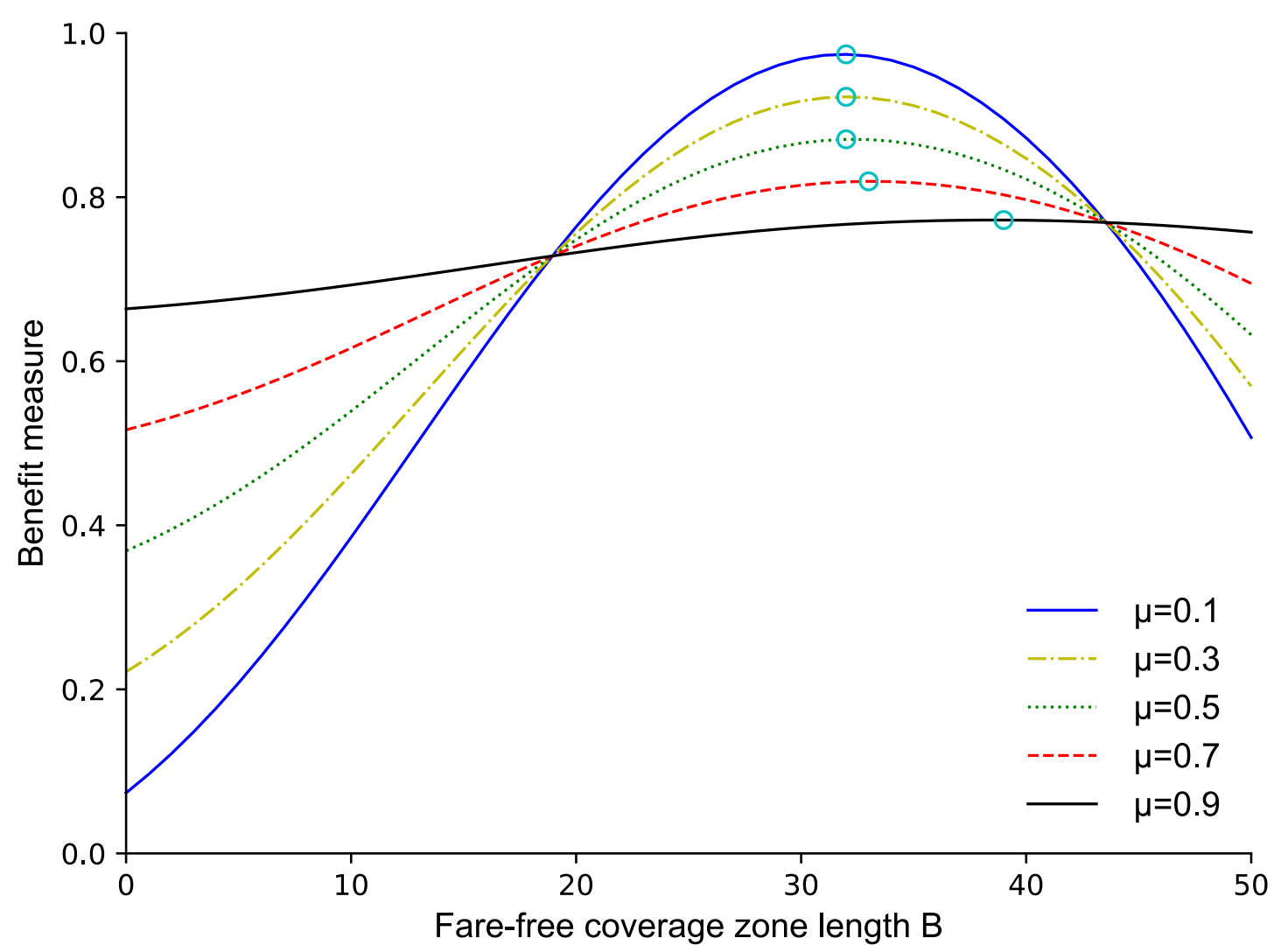

Figure 12: Benefit measure under different weight $\mu$.

## 4.3 Dynamic model analysis

In this section, the proposed dynamic model is applied to a real-world case study in Salt Lake City, Utah. The UTA implemented a fare-free zone in downtown Salt Lake City, allowing passengers to use UTA buses, the TRAX light rail, and later the S-Line streetcar without paying fares (Regional Zero-Fare Study, 2022). However, UTA has not provided $F_F$-service on a complete transit line, though the concept of $F_F$-service has gained increasing attention. Our real-world case study focuses on the 19.3 mile long UTA Blue Light Rail Line (Blue Line), as depicted in Figure 13. This case study offers a practical application of the dynamic model for the Blue Light Rail Line. The line's service frequency has been optimized under conditions of demand uncertainty. The analysis also reveals changes in equity, as measured by the Gini index, for demand groups along the corridor due to the implementation of the fare-free policy.

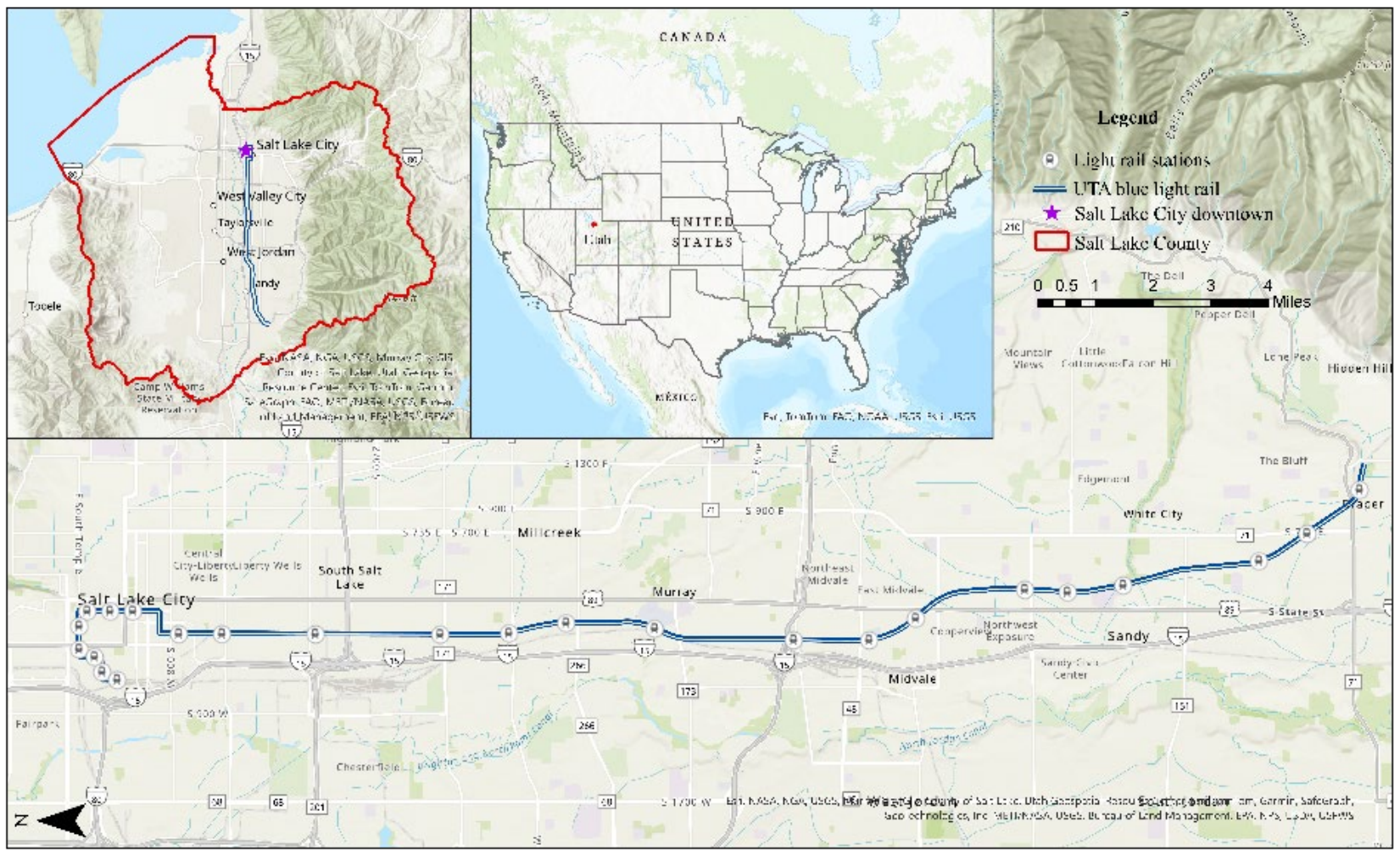

Figure 13: UTA blue light rail corridor in Salt Lake City

The simulation of demand density over time in the CBD is presented by following a GBM model. Ridership data from the UTA Blue Light Rail Line between January 2022 and August 2024 (UTA Ridership Dashboard, 2024) are used to calibrate an annual growth rate $\eta$ of 1.16 % and a volatility rate $\sigma$ of 13.47%.

Based on the length of the Blue Line and the initial mode split, the demand density in the CBD is derived for the Blue Line corridor and is presented in Figure 14, showing historical data trends. Utilizing the estimated annual growth and volatility rates from the historical data, fifty simulated trajectories or paths of the CBD's demand density from September 2024 and February 2028 are generated using Eq. (32), with the simulated data displayed in Figure 14. These fifty trajectories will be used to evaluate the performance of the proposed dynamic policy below.

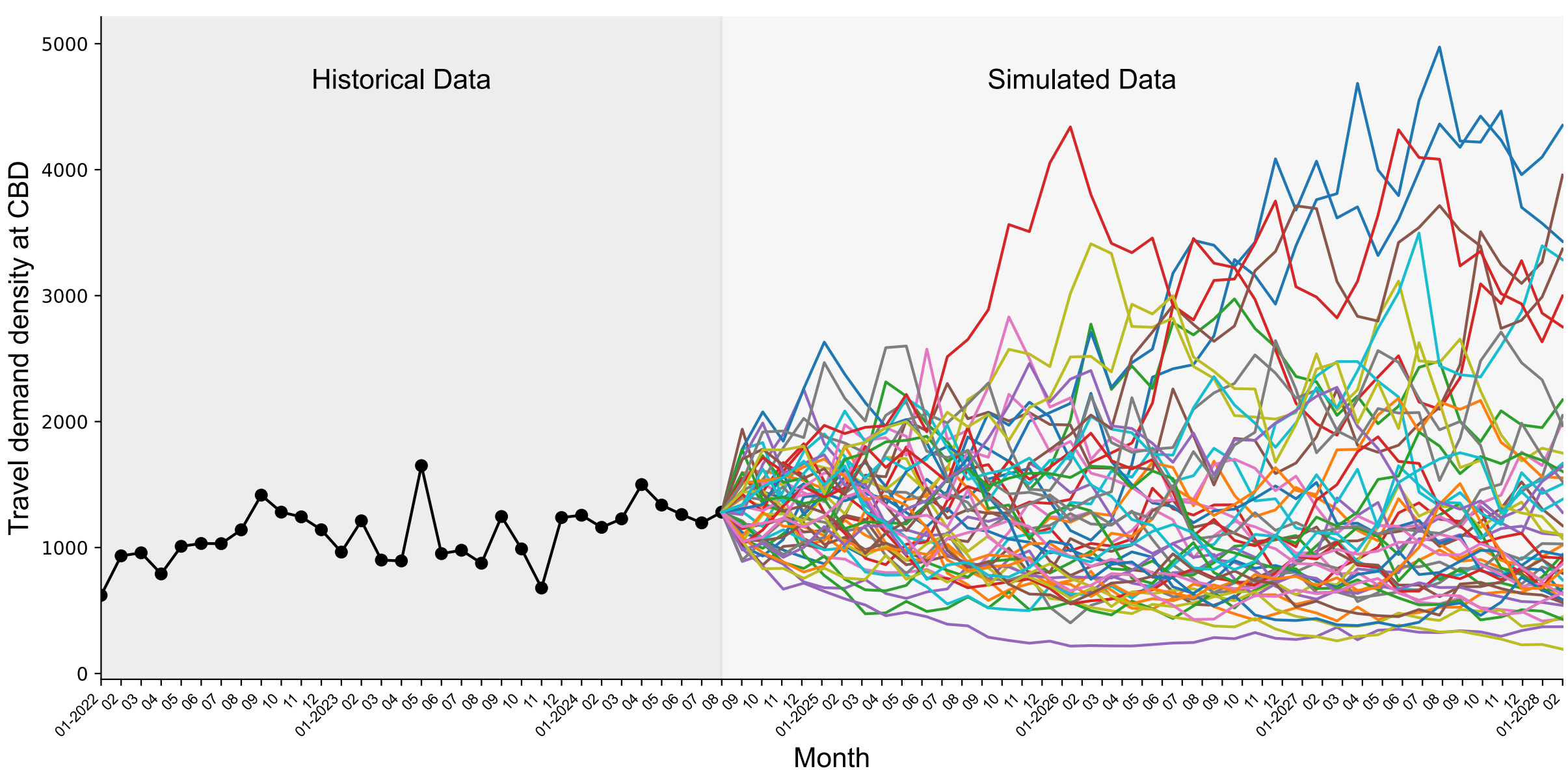

Figure 14: Fifty samples of the travel demand on UTA Blue Line corridor.

Two regimes are analyzed to illustrate how thresholds vary under different benefit-focused objectives: the equity-aware regime (E-regime) and the social welfare-aware regime (W-regime). Using a static model, the most beneficial lengths of $F_F$-service under these regimes are computed and presented in Figure 15. When the operator prioritizes social welfare maximization, the optimal length of $F_F$-service is 16.2, marked by a square. Conversely, when the operator focuses on maximizing social equity, the entire length of the UTA Blue Line should be implemented as $F_F$-service, indicated by a circle.

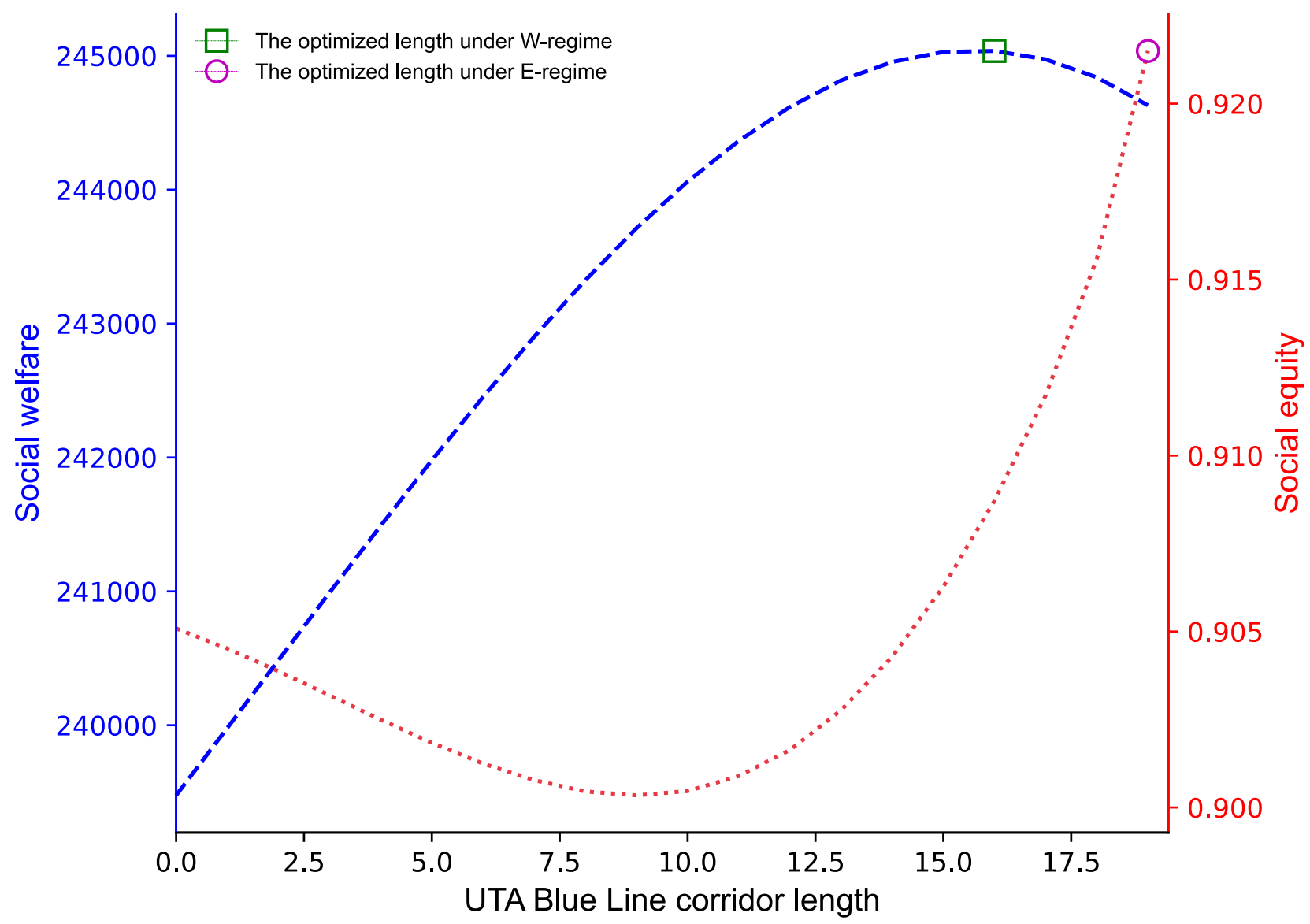

Figure 15. The optimized length of $F_F$-service under different regimes

To obtain $\overline{Q}$ and $\underline{Q}$, the option value $V_0$ and $V_1$ is plotted in Figure 18, which can be computed by Eqs. (41) and (42). The upper threshold $\overline{Q}$ is determined using Eq. (36) when the $V_1(Q)$ exceeds $V_0(Q)$ by switching cost $D$. Similarly, the lower threshold $\underline{Q}$ is derived from Eq. (37) when $V_1(Q)$ falls below $V_0(Q)$ by switching cost $K$. To illustrate the difference between stochastic policies in two cases, Figure 16 and Figure 17 display the thresholds associated with the different stochastic policies (with and without switching cost) for the one simulated trajectory of fifty samples in Figure 15.

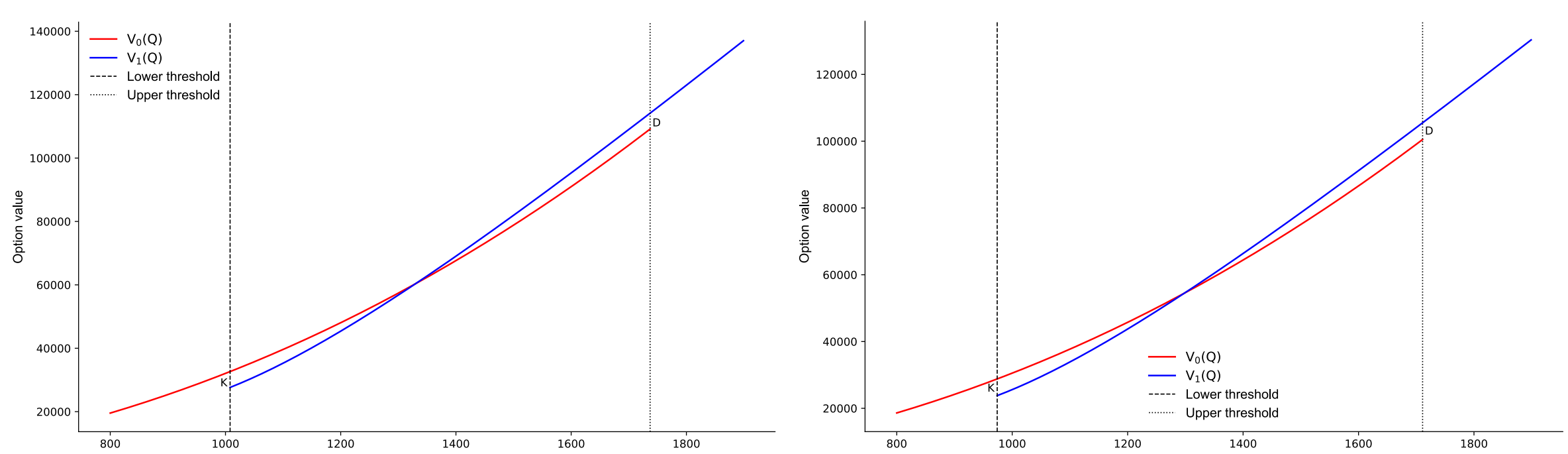

Figure 16: The option value of fare free transit and fare-based transit with respect to demand (a) E-regime policy; (b) W-regime policy.

#### 4.3.1 Social equity-aware regime

As shown in Figure 15, social equity is maximized when free-fare transit is implemented along the entire corridor. Therefore, $F_F$-service should be launched along the entire line if the objective is to maximize social equity. The three thresholds in this case are shown in Figure 17.

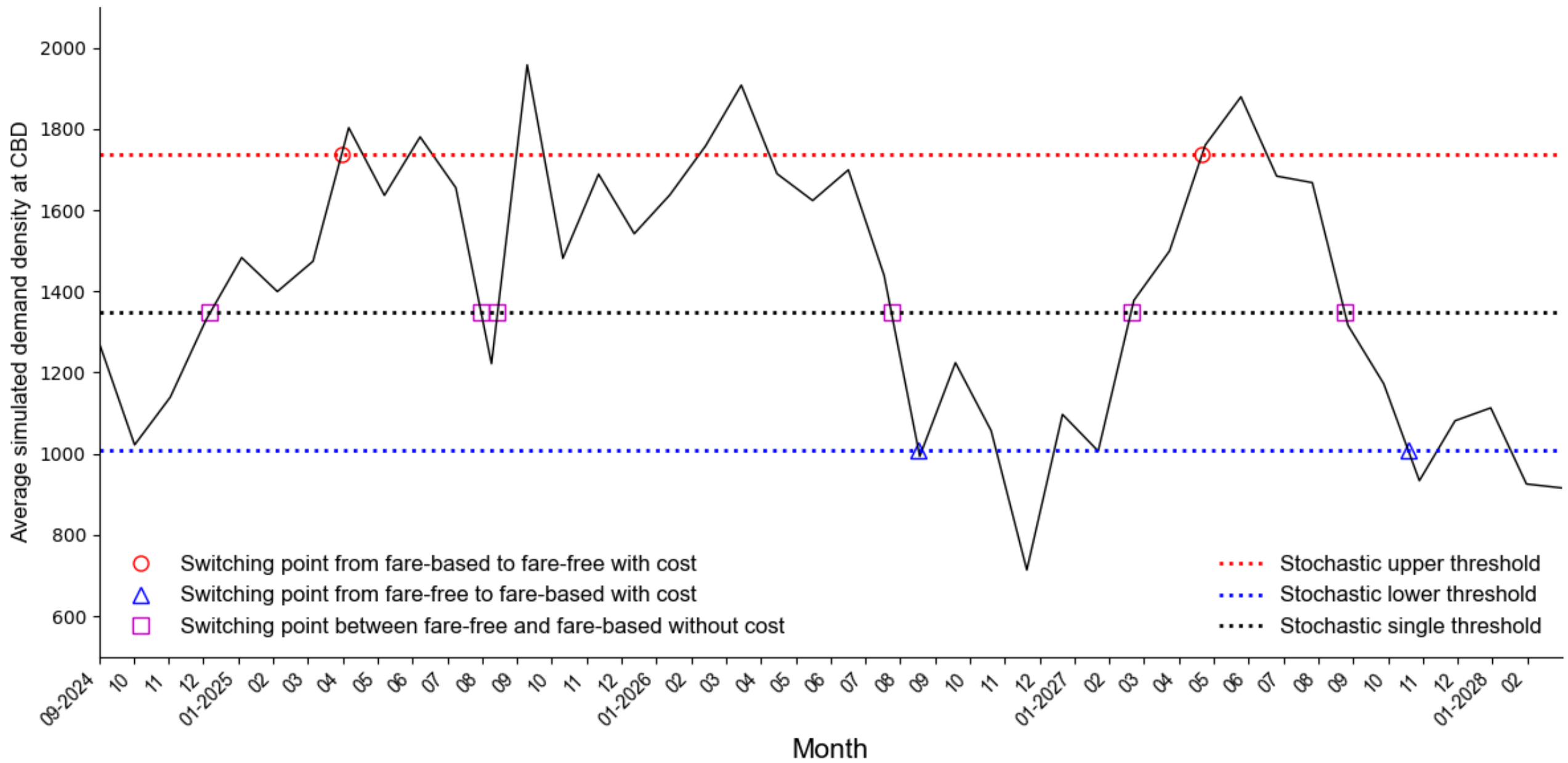

Figure17: Optimized $F_F$-service activation and cancellation timing under E-regime regime.

The optimized upper and lower thresholds under the infinite horizon setting with the estimated stochastic process are determined by Eqs. (45)- (46). In Figure 17, the entry threshold $\overline{Q} = 1737$, represents the demand density where $F_P$-service switches to $F_F$-service. Conversely, the exit threshold, denoted as $\underline{Q} = 1008$, is the demand density at which the $F_F$-service reverts to $F_P$-service. In reverse, the exit threshold is the demand density that switching from $F_F$-service switches to $F_P$-service. If the switching cost is not considered in stochastic model, the one threshold, $Q^* = 1347$ can be obtained by Eq. (47). case. The stochastic case with switching cost exhibits four switching points and the stochastic model without switching cost shows six switching points, represented by squares in Figure 17.

### 4.3.2 Social welfare-aware regime

As shown in Figures 15 in static analysis, the optimized length of $F_P$-service will make the social welfare maximum. Under W-regime, different thresholds are shown in Figure 18. The stochastic switching with a switching cost has four switching points, while the stochastic model without a switching cost displays nine switching points.

Under W-regime, the entry threshold $\overline{Q} = 1711$ indicates the demand density at which the transit service transitions from fare-based to fare-free transit. Similarly, the exit threshold $\underline{Q} = 974$ represents the demand density at which the service switches back from $F_F$-service to $F_P$-service. Conversely, the exit threshold marks the point where the system reverts from $F_F$-service to $F_P$-service. When switching costs are excluded from the stochastic model, a single threshold $Q^* = 1317$ is derived using Eq. (47). Under W-regime, the model with switching costs shows four switching points, whereas the model without switching costs displays six switching points.

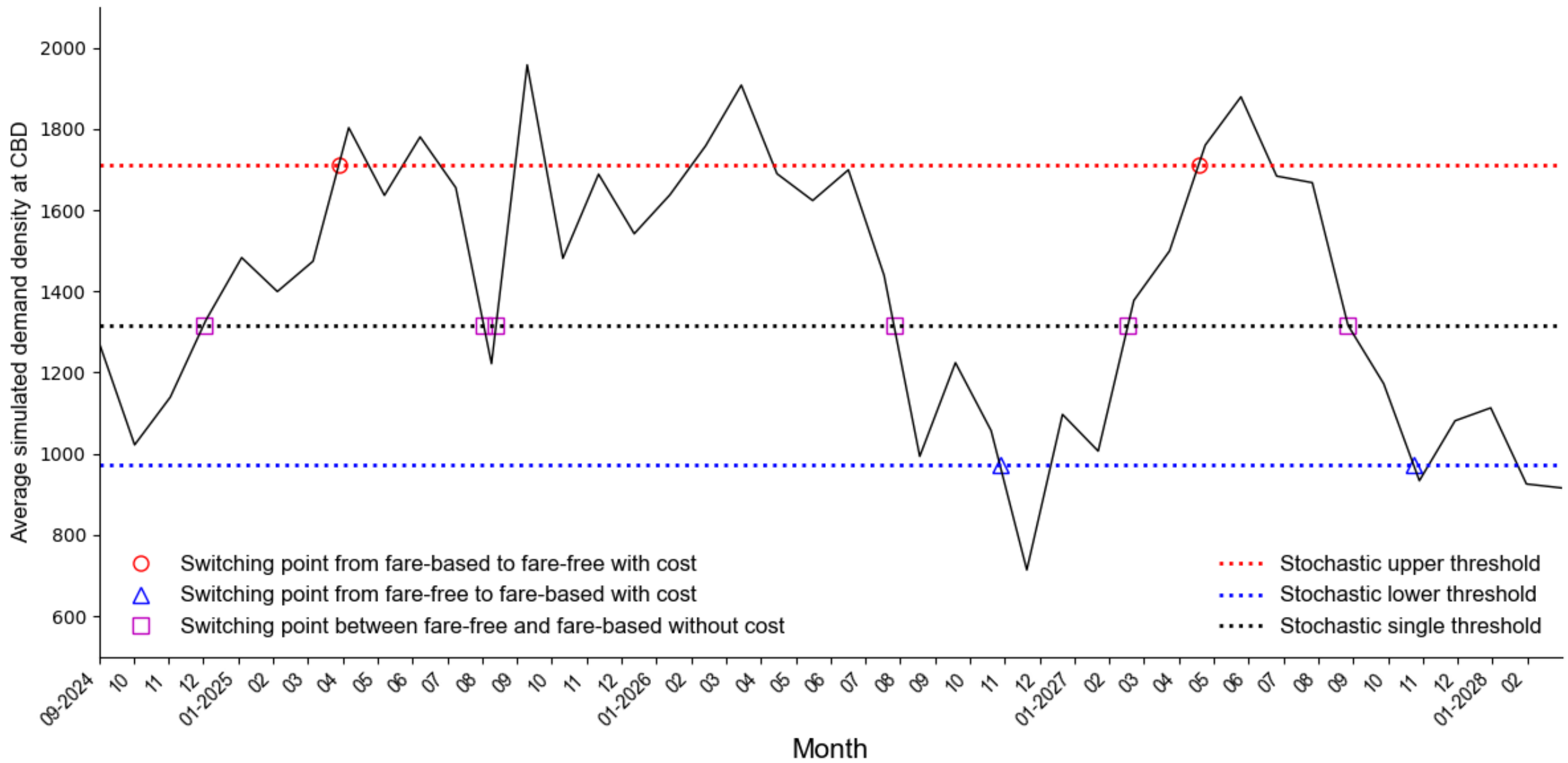

Figure18: Optimized fare-free transit activation and cancellation timing under W-regime

Table 3 compares the thresholds and implementation periods for $F_P$-service and $F_F$-service under E-regime and W-regime. For the W-regime, which optimizes for a 16-mile transit length, the threshold values for $\overline{Q}$ and $\underline{Q}$ are 1711 and 974, respectively. In contrast, the E-regime, which implements $F_F$-service along the full length of UTA Blue Line corridor, shows slightly higher values for $\overline{Q}$ and $\underline{Q}$, with 1737 and 1008, respectively. The changes in these values indicate that lower upper (lower) threshold to active (and cancel) the $F_F$-service are needed under W-regime, which means earlier entry from $F_P$-service to $F_F$-service and later exit from $F_F$-service to $F_P$-service, as shown in Figure 17 and 18. This is because the social welfare of implementing $F_F$-service increases under W-regime (optimized length) compared to E-regime (full length), which can be shown in Figure 15 and Table 3. In terms of operational planning over the horizon, implementing $F_F$-service is beneficial during two specific periods, while in three other periods, it is preferable to maintain a $F_P$-service. Although the number of implementing periods for $F_F$-service and $F_P$-service remains unchanged between E-regime and W-regime, the total duration of $F_F$-service is extended from 22 months under E-regime to 26 months under W-regime, as shown in the implementation periods of the $F_F$-service in Table 3. The threshold decreases from 1,347 in the E-regime to 1,317 in the W-regime when switching costs are excluded. By strategically alternating between these models based on demand density, the service can optimize operational efficiency and better match demand fluctuations throughout the planning horizon.

Furthermore, the average values of social welfare and implantation duration across all samples are computed by applying the obtained switching thresholds to each of the 50 samples under both regimes. The results show that the W-regime achieves a slightly higher average social welfare than the E-regime, which highlights the social welfare improvement achieved through optimized service length. Additionally, the average implementation duration of the $F_P$-service decreases from 21.4 months under the E-regime to 17.8 months under the W-regime, whereas the average duration of the $F_F$-service increases from 20.6 to 24.2 months.

Table 3: The difference of thresholds and implementation timing under different regimes.

| | Fare-free transit length | Switching threshold | | One sample case | | | Average value across fifty samples | | |
|---|---|---|---|---|---|---|---|---|---|
| | $B$ | $Q^*$ | $\overline{Q}$ | $\underline{Q}$ | $W$ | Implementation periods in planning time: $F_P$-service | Implementation periods in planning time: $F_F$-service | $\overline{W}$ | $F_P$-service | $F_F$-service |
| E-regime | Full length (19.4 miles) | 1347 | 1737 | 1008 | 244723 | 09/2024-04/2025; 08/2026-05/2027; 10/2027-02/2028. (3 periods (Total 20 months)) | 04/2025-08/2026; 05/2027-10/2027. (2 periods (Total 22 months)) | 243775 | 21.4 months | 20.6 months |
| W-regime | Optimized length (16 miles) | 1317 | 1711 | 974 | 245472 | 09/2024-03/2025; 11/2026-04/2027; 10/2027-02//2028. (3 periods (Total 16 months)) | 03/2025-11/2026; 04/2027-10/2027. (2 periods (Total 26 months)) | 245328 | 17.8 months | 24.2 months |

# 5. Conclusion and Future Studies

In this study, an optimized scheme for fare-free transit services is explored, focusing on determining the length of the coverage zone and transit frequency, and establishing the timing for starting and stopping the policy. A market entry and exit model is proposed, grounded on the principles of geometric Brownian motion, to tackle the timing of the implementation and cessation of a fare-free transit service, incorporating the inherent uncertainty in demand density. The optimized zone length is derived numerically, whereas the timing of fare-free activation and cessation is examined analytically within the contexts of fare-based and fare-free transit services. The social equity among different demand groups is quantified. Notably, under a social equity consideration, the optimized zone length may diverge from the length that maximizes social welfare.

The static and dynamic models yield significant insights and quantitative findings:

Firstly, establishing a fare-free zone significantly promotes bus demand, while outside the zone, the bus share decreases and auto share increases.

Secondly, as the fare-free zone expands, both bus and auto users within the zone increase, whereas those outside the zone decrease.

Thirdly, activating fare-free transit services boosts the social equity ratio, with benefits increasing as transit fares decrease. Based on the Lorenz curve, fare-free transit in the whole corridor yields the highest level of equity compared with the fare-free zone and fare-based transit.

Fourthly, the benefit measure shows that when considering both social equity and social welfare, the optimized zone length increases compared to focusing solely on social welfare. Specifically, incorporating social equity into the benefit measure requires sacrificing some social welfare benefits, as implementing fare-free transit along the entire corridor achieves social equity maximization.

Finally, the most beneficial fare-free transit length for the UTA Blue Line under different regimes is identified using static models. When applying dynamic models, the implementation timing can also be determined by accounting for demand uncertainty. Compared to the E-regime, the three thresholds decrease under W-regime. This results in earlier implementation and longer durations of fare-free transit under W-regime. In scenarios where the switching cost between the two services is negligible, a single threshold emerges. Without switching cost, the frequency of switching increases since such costs are not a consideration. Under E-regime, the switching times are higher compared to the W-regime.

The proposed modeling methodology can provide a valuable tool for decision-making regarding the operational scheme of $F_F$-service, when maximizing social welfare while also considering social equity. The derived closed-form solution for implementation length and timing under different regimes contributes valuable insights to the discussion on $F_F$-service.

Future research endeavors in the sphere of fare-free transit time periods and spatially distributed zones may explore the following.

Primarily, the contemporary post-COVID era has drastically altered travel patterns, necessitating an understanding of peak transit demand times and locations. Future studies may consider a many-to-many travel pattern, as many persons still engage in remote work. Equity, an integral concern in public transit, mandates an examination of fare-free periods and zones, ensuring they adequately serve disadvantaged and marginalized communities. In the quantification of equity, it would be desirable to account for diverse household distribution brackets.

Moreover, a comprehensive analysis of the financial implications of fare-free periods/zones is desirable. This involves not merely accounting for the loss of fare revenue, but also potential increases in ridership and reduced costs of fare collection. It would be desirable to concurrently analyze transit operation parameters, such as frequency, fare, and vehicle size, when optimizing the activation timing and zoning of fare-free.

Future investigations should also incorporate environmental benefits and impacts on traffic congestion. If $F_F$-service can appreciably curtail car usage, especially during peak periods, the environmental and societal benefits may potentially outweigh the associated costs.

Consideration of congestion effects could be beneficial given the interaction of multiple modes of transport. Future research should focus on the interplay between fare-free periods/zones and other modes of transport such as bicycling and car sharing, to conceptualize a cohesive, efficient, and equitable transportation system.

In quantifying the uncertainty of fare collection cost, the development of novel statistical models and techniques would be beneficial. Stochastic models, in particular, can incorporate multiple uncertainties, e.g., regarding fare collection cost, operating cost, and maintenance cost, thus leading to a more realistic model.

## Appendix A

Note that $\lim_{Q\to 0+} V_0(Q) = 0$; thus we must have $X_0 = 0$. Similarly, considering $\lim_{Q\to\infty} V_1(Q) = 0$, we have $Y_1 = 0$ (Dixit,1994). Then we can omit the subscripts on the remaining coefficients and write the solutions as Eqs. (43) and (44):

$$V_0(Q) = Y_0 Q^{\gamma_1}, \tag{A2}$$

and

$$V_1(Q) = X_1 Q^{\gamma_0} + \frac{\left(d_1^Q - d_0^Q\right)}{k-\eta} Q - \frac{\left(d_1^C - d_0^C\right)}{k}, \tag{A3}$$

For specific choices of $\eta$, $\sigma$, and $k$, Eqs. (A2) and (A3) can be substituted into the value matching and smooth pasting conditions in Eq. (37)-(40) to establish four equations and four unknown variables: $Y_0$, $X_1$, $\overline{Q}$, and $\underline{Q}$. It can be shown that the optimized solution $\left[Y_0, X_1, \overline{Q}, \underline{Q}\right]$ is uniquely determined by solving the system of nonlinear equations in (A4):

$$\begin{bmatrix} -Y_0\overline{Q}^{\gamma_1} + X_1\overline{Q}^{\gamma_0} + \frac{\left(d_1^Q - d_0^Q\right)}{k-\eta}\overline{Q} - \frac{\left(d_1^C - d_0^C\right)}{k} - D \\ -Y_0\underline{Q}^{\gamma_1} + X_1\underline{Q}^{\gamma_0} + \frac{\left(d_1^Q - d_0^Q\right)}{k-\eta}\underline{Q} - \frac{\left(d_1^C - d_0^C\right)}{k} + K \\ -Y_0\gamma_1\overline{Q}^{\gamma_1-1} + X_1\gamma_0\overline{Q}^{\gamma_0-1} + \frac{\left(d_1^Q - d_0^Q\right)}{k-\eta} \\ -Y_0\gamma_1\underline{Q}^{\gamma_1-1} + X_1\gamma_0\underline{Q}^{\gamma_0-1} + \frac{\left(d_1^Q - d_0^Q\right)}{k-\eta} \end{bmatrix} = 0. \tag{A4}$$

After rearranging Eq. (A4), the upper threshold and lower threshold can be obtained with Eqs. (43) and (44).

## Appendix B

For the special case where the switching cost becomes zero ($D = K = 0$), the upper and lower thresholds $\overline{Q}$ and $\underline{Q}$ converge to a single threshold $Q^*$. Then Eqs. (37)-(38) become

$$V_0(Q^*) = V_1(Q^*), \tag{B1}$$

$$V_0'(Q^*) = V_1'(Q^*). \tag{B2}$$

That is,

$$\begin{bmatrix} Y_0 Q^{*\gamma_1} - X_1 Q^{*\gamma_0} - \nu_0 Q^* + \nu_1 \\ Y_0 \gamma_1 Q^{*\gamma_1 - 1} - X_1 \gamma_0 Q^{*\gamma_0 - 1} - \nu_0 \end{bmatrix} = 0. \tag{B3}$$

Note that we have to determine three unknown parameters. To find $Q^*$, we also need the following condition

$$V_0''(Q^*) = V_1''(Q^*). \tag{B4}$$

From Eq. (B4), we have

$$Y_0 \gamma_1 (\gamma_1 - 1) Q^{*(\gamma_1 - 2)} = X_1 \gamma_0 (\gamma_0 - 1) Q^{*(\gamma_0 - 2)}. \tag{B5}$$

Solving these non-linear equations (B1), (B2) and (B5), we obtain $Q^*$ as Eq. (45).

## References

Balliauw, M., Kort, P. M., Meersman, H., Van de Voorde, E., & Vanelslander, T. (2020). The case of public and private ports with two actors: capacity investment decisions under congestion and uncertainty. *Case Studies on Transport Policy*, 8(2), 403-415.

Ben-Elia, E., & Benenson, I. (2019). A spatially-explicit method for analyzing the equity of transit commuters' accessibility. *Transportation Research Part A: Policy and Practice*, *120*, 31-42.

Boyd, B., Chow, M., Johnson, R., & Smith, A. (2003). Analysis of effects of fare-free transit program on student commuting mode shares: BruinGo at University of California at Los Angeles. *Transportation Research Record*, *1835*(1), 101-110.

Brennan, M. J., & Schwartz, E. S. (1978). Finite difference methods and jump processes arising in the pricing of contingent claims: A synthesis. *Journal of Financial and Quantitative Analysis*, *13*(3), 461-474.

Bull, O., Muñoz, J. C., & Silva, H. E. (2021). The impact of fare-free public transport on travel behavior: Evidence from a randomized controlled trial. *Regional Science and Urban Economics*, 86, 103616.

Butler, A., & Sweet, M. (2020). No free rides: Winners and losers of the proposed Toronto Transit Commission U-Pass program. *Transport Policy*, *96*, 15-28.

Button, K. (2010). *Transport economics*. Edward Elgar Publishing.

Cantillo, A., Raveau, S., & Muñoz, J. C. (2022). Fare evasion on public transport: Who, when, where and how? *Transportation Research Part A: Policy and Practice*, *156*, 285-295.

Cats, O., Susilo, Y. O., & Reimal, T. (2017). The prospects of fare-free public transport: evidence from Tallinn. *Transportation*, 44(5), 1083-1104.

Cats, O., Susilo, Y., & Eliasson, J. (2012). Impacts of Free PT, Tallinn–Evaluation Framework. *Department of Transport Science in the Royal Institute of Technology (KTH),* Sweden. (October). http://www. tallinn. ee/eng/tasutauhistransport/g9616s62872.

Chang, S. K., & Schonfeld, P. M. (1993). Welfare maximization with financial constraints for bus transit systems. *Transportation Research Record*, (1395).

Chen, R., & Zhou, J. (2022). Fare adjustment's impacts on travel patterns and farebox revenue: An empirical study based on longitudinal smartcard data. *Transportation Research Part A: Policy and Practice*, *164*, 111-133.

Chen, Z., Guo, Y., Stuart, A. L., Zhang, Y., & Li, X. (2019). Exploring the equity performance of bike-sharing systems with disaggregated data: A story of southern Tampa. *Transportation research part A: policy and practice*, *130*, 529-545.

Chow, J. Y. J. (2016). Dynamic UAV-based traffic monitoring under uncertainty as a stochastic arc-inventory routing policy. *International Journal of transportation science and technology*, *5*(3), 167-185.

Chow, J. Y. J. (2018). *Informed Urban transport systems: Classic and emerging mobility methods toward smart cities*. Elsevier.

Chow, J. Y. J., & Regan, A. C. (2011a). Network-based real option models. *Transportation Research Part B: Methodological*, 45(4), 682-695.

Chow, J. Y. J., & Regan, A. C. (2011b). Real option pricing of network design investments. *Transportation Science*, 45(1), 50-63.

Chow, J. Y. J., & Sayarshad, H. R. (2016). Reference policies for non-myopic sequential network design and timing problems. *Networks and Spatial Economics*, *16*, 1183-1209.

Chow, J. Y. J., Regan, A. C., Ranaiefar, F., & Arkhipov, D. I. (2011c). A network option portfolio management framework for adaptive transportation planning. *Transportation Research Part A: Policy and Practice*, 45(8), 765-778.

Couto, G., Nunes, C., & Pimentel, P. (2015). High-speed rail transport valuation and conjecture shocks. *The European Journal of Finance*, 21(10-11), 791-805.

Cox, J. C., Ross, S. A., & Rubinstein, M. (1979). Option pricing: A simplified approach. *Journal of financial Economics*, *7*(3), 229-263.

Dai, J., Liu, Z., & Li, R. (2021). Improving the subway attraction for the post-COVID-19 era: The role of fare-free public transport policy. *Transport Policy*, 103, 21-30.

De Witte, A., Macharis, C., & Mairesse, O. (2008). How persuasive is 'free'public transport?: a survey among commuters in the Brussels Capital Region. *Transport Policy*, *15*(4), 216-224.

De Witte, A., Macharis, C., Lannoy, P., Polain, C., Steenberghen, T., & Van de Walle, S. (2006). The impact of "free" public transport: The case of Brussels. *Transportation Research Part A: Policy and Practice*, *40*(8), 671-689.Delbosc, A., & Currie, G. (2011). Using Lorenz curves to assess public transport equity. *Journal of Transport Geography*, *19*(6), 1252-1259.

Dixit, A. (1989). Entry and exit decisions under uncertainty. *Journal of political Economy*, *97*(3), 620-638.

Dixit, R. K., Dixit, A. K., & Pindyck, R. S. (1994). *Investment under uncertainty*. Princeton University Press.

Errami, M., Russo, F., & Vallois, P. (2002). Itô's formula for C 1, λ-functions of a càdlàg process and related calculus. *Probability theory and related fields*, *122*, 191-221.

Farber, S., Bartholomew, K., Li, X., Páez, A., & Habib, K. M. N. (2014). Assessing social equity in distance based transit fares using a model of travel behavior. *Transportation Research Part A: Policy and Practice*, *67*, 291-303.

Fujita, M. (1989). *Urban economic theory*. Cambridge Books.

Gao, J., & Li, S. (2024). Synergizing shared micromobility and public transit towards an equitable multimodal transportation network. *Transportation Research Part A: Policy and Practice*, *189*, 104225.

Giversen, J., & Bendkia, M. (2011). Continuous Time Processes in Times of Crisis: The Case Of GBM and CEV Models.

Goldberg, D. (2021). The New Geography of Fare-Free Public Transport: Examining the Suspension of Fares in the United States During COVID-19. Brussels: Master in Urban Studies (ULB/VUB).

Gray, A., 2018. Estonia is Making Public Transport Free. Available from: https://www.weforum.org/agenda/2018/06/estonia-is-making-public-transport-free/

Guo, Q. W., Chow, J. Y. J., & Schonfeld, P. (2018). Stochastic dynamic switching in fixed and flexible transit services as market entry-exit real options. *Transportation Research Part C: Emerging Technologies*, *94*, 288-306.

Guo, Q., Chen, S., Sun, Y., & Schonfeld, P. (2023). Investment timing and length choice for a rail transit line under demand uncertainty. *Transportation Research Part B: Methodological*, *175*, 102800.

Guo, Q., Sun, Y., Schonfeld, P., & Li, Z. (2021). Time-dependent transit fare optimization with elastic and spatially distributed demand. *Transportation Research Part A: Policy and Practice*, 148, 353-378.

Hodge, D. C., Orrell III, J. D., & Strauss, T. R. (1994). Fare-free policy: costs, impacts on transit service, and attainment of transit system goals. Washington State Department of Transportation Report. http://www.wsdot.wa.gov/research/reports/fullreports/277.1.pdf

Kębłowski, W. (2020). Why (not) abolish fares? Exploring the global geography of fare-free public transport. *Transportation*, *47*(6), 2807-2835.

Kirschen, M., Pettine, A., Adams, M., Persaud, H., 2022. Fare-free transit evaluation framework. Transportation Research Board. https://doi.org/10.17226/26732.

Li, Y., & DaCosta, M. N. (2013). Transportation and income inequality in China: 1978–2007. *Transportation Research Part A: Policy and Practice*, *55*, 56-71.

Li, Z. C., Guo, Q. W., Lam, W. H., & Wong, S. C. (2015). Transit technology investment and selection under urban population volatility: A real option perspective. *Transportation Research Part B: Methodological*, 78, 318-340.

Longstaff, F. A., & Schwartz, E. S. (2001). Valuing American options by simulation: a simple least-squares approach. *The review of financial studies*, *14*(1), 113-147.

Lu, Q. L., Mahajan, V., Lyu, C., & Antoniou, C. (2024). Analyzing the impact of fare-free public transport policies on crowding patterns at stations using crowdsensing data. *Transportation Research Part A: Policy and Practice*, *179*, 103944.

Noël, T. S. F. (2022). *Fare-free public transportation effects' assessment on french economy* (Doctoral dissertation).

Powell, W. B. (2007). *Approximate Dynamic Programming: Solving the curses of dimensionality* (Vol. 703). John Wiley & Sons.Rath, S., & Chow, J. Y. J. (2022). A deep real options policy for sequential service region design and timing. *arXiv preprint arXiv:2212.14800*.

Regional Zero-Fare Study Final Report: https://wfrc.org/Studies/ZeroFareTransit/ZeroFareTransitStudy_FinalReport.pdf

Saphores, J.D., Shah D., and Khatun, F. (2020). A review of reduced and free transit fare programs in California. The University of California Institute of Transportation Studies Report. https://doi.org/10.7922/G2XP735.

Sødal, S., Koekebakker, S., & Aadland, R. (2008). Market switching in shipping—A real option model applied to the valuation of combination carriers. *Review of Financial Economics*, *17*(3), 183-203.

Studenmund, A. H., & Connor, D. (1982). The free-fare transit experiments. *Transportation Research Part A: General*, *16*(4), 261-269.

Tirachini, A., & Hensher, D. A. (2011). Bus congestion, optimal infrastructure investment and the choice of a fare collection system in dedicated bus corridors. *Transportation Research Part B: Methodological*, 45(5), 828-844.

UTA Ridership Dashboard 2024 UTA Mode Level Boardings Weekday Averages | UTA Mode Level Boardings Weekday Averages | UTA Open Data Portal (rideuta.com)

Verheyen, G., 2010. Irrationality in Modal Choices: Free Public Transport. Master's thesis. Hasselt University. Diepenbeek, Belgium. Original dutch title: ''Irrationaliteit in Vervoerswijzekeuze: Gratis".

Volinski, J. (2012). *Implementation and outcomes of fare-free transit systems* (No. 101). Transportation Research Board.

Wang, Q., & Deng, L. (2019). Integrated optimization method of operational subsidy with fare for urban rail transit. *Computers & Industrial Engineering*., 127, 1153–1163.

Wang, Q., Schonfeld, P., Deng, L., Xu, G., & Ling, S. (2023). Optimization of differentiated fares and subsidies for different urban rail transit users. *Computers & Industrial Engineering*, 179, 109144.

Wang, S., Liu, Y., & Corcoran, J. (2021). Equity of public transport costs before and after a fare policy reform: An empirical evaluation using smartcard data. *Transportation Research Part A: Policy and Practice*, *144*, 104-118.

Zhang, F., Yang, Z., Jiao, J., Liu, W., & Wu, W. (2020). The effects of high-speed rail development on regional equity in China. *Transportation Research Part A: Policy and Practice*, *141*, 180-202.

Zhou, J., Zhang, M., & Zhu, P. (2019). The equity and spatial implications of transit fare. *Transportation Research Part A: Policy and Practice*, *121*, 309-324.